\documentclass{article}
\usepackage{etex}
\usepackage{color}
\usepackage[latin1]{inputenc}
\usepackage{amsmath, enumerate, amsfonts}
\usepackage[widemargins]{a4}
\usepackage{amsmath, enumerate, amsfonts,amssymb, pictexwd, amscd}

\newcommand{\pcom}{^\wedge_p}
\usepackage[all]{xy}
\xyoption{arc}
\newtheorem{Conjecture}{Conjecture}[section]

\newtheorem{Theorem}[Conjecture]{Theorem}
\newtheorem{Definition}[Conjecture]{Definition}
\newtheorem{Proposition}[Conjecture]{Proposition}

\newtheorem{Remark}[Conjecture]{Remark}

\newcommand*{\F}{\mathbb{F}}

\newcommand*{\hoTop}{\textbf{hoTop}}
\DeclareMathOperator{\Aut}{Aut}

\DeclareMathOperator{\Mor}{Mor}
\DeclareMathOperator{\hocolim}{hocolim}
\usepackage[latin1]{inputenc}
\usepackage{graphicx}
\usepackage{amsmath, enumerate, amsfonts}
\usepackage[widemargins]{a4}
\usepackage[all]{xy}
\xyoption{arc}

\title{Assigning a classifying space to a fusion system up to F-isomorphism}
\author{Nora Seeliger}
\date{}
\begin{document}
\maketitle
\begin{abstract}
Complementing and extending the Inventiones work of Benson, Grodal, Henke [Group cohomology and control of $p-$fusion, Invent.\ Math. 197 (2014), 491--507] we give criteria for a space to have cohomology (strongly) $F-$isomorphic in the sense of Quillen to the stable elements. We extend results about groups models to fusion systems over discrete $p$-toral groups and profinite groups and provide various applications to the Kan-Thurston Theorem.
\end{abstract}
\section{Introduction}
\label{cIntro}
In the topological theory of $p-$local finite groups introduced by Broto, Levi, and Oliver one tries to approximate the classifying space of a finite group via the $p-$local structure of the group, at least up to $\mathbb{F}_p-$cohomology. Every finite group gives canonically rise to a saturated fusion system for every prime dividing its order however not every fusion system is the fusion system of a finite group. In contrast there exist infinite groups realizing arbitrary fusion systems. In this article we give an abstract criterion for a space to have cohomology $F-$isomorphic and strongly $F-$isomorphic in the sense of Quillen to the stable elements. Classifying spaces of certain group models are examples of such spaces. This is one motivation for our continued study of group models and in this context we give an illustration of Quillen's theorem. We finish with extending classical results about group models to fusion systems over discrete $p-$toral groups and pro-fusion systems over pro $p-$groups.\\
The author was supported by
ANR BLAN08-2-338236, an
Erwin-Schr\"odinger-Institute Junior-Research-Fellowship, a
Mathematisches-Forschungsinstitut-Oberwolfach Leibniz-Fellowship, an
invitation to the Max-Planck-Institut f\"ur Mathematik in Bonn, ARC Discovery Project DP120101399 at the Australian National University and a postdoctoral fellowship at the University of Haifa.
I would like to thank David
Benson, David Blanc, David Green, Hans-Werner Henn,
 Ian Leary,
Ran Levi, Ehud Meir, Sune Precht Reeh, and Tomer Schlank for enriching discussions on the topic, and in addition David Blanc and Ehud Meir for giving comments on earlier versions of the manuscript.
\section{Preliminaries}
\subsection{Fusion Systems}
We review basic definitions of fusion systems and
linking systems and establish our notations. Our main references are
\cite{BLO2}, \cite{BLO4} and \cite{IntroMarkus}. Let $S$ be a finite
$p$-group. A \textbf{fusion system} $\mathcal{F}$ over $S$ is a
category whose objects are all the subgroups of $S$, and which
satisfies the following two properties for all $P,Q\leq S$: The set
$Hom_{\mathcal{F}}(P,Q)$ contains injective group homomorphisms and
amongst them all morphisms induced by conjugation of elements in $S$
and each element is the composite of an isomorphism in $\mathcal{F}$
followed by an inclusion. Two subgroups $P,Q\leq S$ will be called
\textbf{$\mathcal{F}-$conjugate} if they are isomorphic in
$\mathcal{F}$. Define
$Out_{\mathcal{F}}(P)=Aut_{\mathcal{F}}(P)/Inn(P)$ for all $P\leq
S$.
 A subgroup $P\leq S$ is \textbf{fully centralized} resp.\ \textbf{fully normalized} in $\mathcal{F}$ if $|C_S(P)|\geq |C_S(P')|$ resp.\ $|N_S(P)|\geq |N_S(P')|$ for all $P'\leq S$ which is $\mathcal{F}$-conjugate to $P$.
 $\mathcal{F}$ is called \textbf{saturated} if for all $P\leq S$ which is fully normalized in $\mathcal{F}$, $P$ is fully centralized in $\mathcal{F}$ and $Aut_S(P)\in Syl_p(Aut_{\mathcal{F}}(P))$ and moreover if $P\leq S$ and $\phi\in Hom_{\mathcal{F}}(P,S)$ are such that $\phi (P)$ is fully centralized, and if we set
$N_{\phi}=\{g\in N_S(P)|\phi c_g \phi^{-1}\in Aut_S(\phi (P))\}$,
then there is $\overline{\phi}\in Hom_{\mathcal{F}}(N_{\phi},S)$
such that $\overline{\phi}|_P=\phi$. A subgroup $P\leq S$ will be
called $\mathcal{F}$-\textbf{centric} if $C_S(P')\leq P'$ for all
$P'$ which are $\mathcal{F}-$conjugate to $P$. Denote
$\mathcal{F}^c$ the full subcategory of $\mathcal{F}$ with objects
the $\mathcal{F}$-centric subgroups of $S$. 
Let
$\mathcal{O}(\mathcal{F})$ be the \textbf{orbit category} of
$\mathcal{F}$ with objects the same objects as $\mathcal{F}$ and morphisms the set
 $Mor_{\mathcal{O}(\mathcal{F})}(P,Q)=Mor_{\mathcal{F}}(P,Q)/Inn(Q)$.
 Let $\mathcal{O}^c(\mathcal{F})$ be the full subcategory of
 $\mathcal{O}(\mathcal{F})$ with objects the $\mathcal{F}-$centric
 subgroups of $\mathcal{F}$. 
 Let $G$ be a discrete group. A finite subgroup $S$ of $G$ will be
called a \textbf{Sylow} $p$-subgroup of $G$ if $S$ is a $p$-subgroup
of $G$ and all $p-$subgroups of $G$ are conjugate to a subgroup of $S$. A discrete group $G$ has \textbf{finite virtual cohomological dimension} if it has a torsion free subgroup of finite index. A group $G$ is called $p$\textbf{-perfect} if it has no normal subgroup of index $p$. Let $\mathcal{F}$ be a saturated fusion
system over the finite $p$-group $S$. A subgroup $P\leq S$ is called
$\mathcal{F}$-\textbf{radical} if
$O_p(Aut_{\mathcal{F}}(P))=Aut_P(P)$, where $O_p(-)$ denotes the
maximal normal $p-$subgroup. A subgroup $P\leq S$ which is both $\mathcal{F}$-centric and $\mathcal{F}$-radical is called $\mathcal{F}$-centric-radical. Denote
$\mathcal{F}^{cr}$ the full subcategory of $\mathcal{F}$ with objects
the $\mathcal{F}$-centric-radical subgroups of $S$. A subgroup $P\leq S$ is called
$\mathcal{F}$\textbf{-essential} if P is $\mathcal{F}$-centric and
$Aut_{\mathcal{F}}(P)/Aut_P(P)$ has a strongly $p-$embedded
subgroup, i.\ e.\ there exists a proper subgroup $Q<P$ which contains
a Sylow $p-$subgroup $R$ of $P$ and $R\neq 1$ but $Q\cap ^xR=1$ for any
$x\in P-Q$. A subgroup $Q\triangleleft S$ is \textbf{normal} in $\mathcal{F}$ if each $\alpha \in Hom_{\mathcal{F}(P,P')}$ extends to a morphism $\overline{\alpha}\in Hom_{\mathcal{F}}(PQ,P'Q)$ which sends $Q$ to itself. A saturated fusion system is called \textbf{constrained} if there is some $Q\triangleleft S$ which is $\mathcal{F}$-centric and normal in $\mathcal{F}$.
Let $S$ be a finite $p-$group and let $P_1,...,P_r,Q_1,...,Q_r$ be
subgroups of $S$. Let $\phi _1,...,\phi _r$ be injective group
homomorphisms $\phi _i:P_i\rightarrow Q_i$ $\forall i$. The fusion
system \textbf{generated} by $\phi _1,...,\phi_r$ is the minimal
fusion system $\mathcal{F}$ over $S$ containing $\phi _1,...,\phi
_r$. Let $P_1,\dots,P_n$ be a collection of $\mathcal{F}$-centric
subgroups of $S$, which might contain repetitions.
Let $K_1,\dots K_n$ be groups with $K_i$ is a subgroup of $\Aut_{\mathcal{L}}(P_i)$
containing $\delta (P_i)$ for all $i=1,\cdots ,n$. 
We say that $K=\{K_1,\dots K_n\}$ is  \textbf{generating} if $\pi(K_1),\dots,\pi(K_n)$ generate
$\mathcal{F}$. Let $\mathcal{F}$ be a fusion system on a finite
$p-$group $S$. A subgroup $T\leq S$ is \textbf{strongly closed} in
$S$ with respect to $\mathcal{F}$, if for each subgroup $P$ of $T$,
each $Q\leq S$, and each $\phi\in Mor_{\mathcal{F}}(P,Q)$, $\phi
(P)\leq T$.
Fix any pair $S\leq G$, where $G$ is a discrete group and $S$ is a finite $p-$subgroup.
    Define $\mathcal{F}_S(G)$ to be the category whose objects are the subgroups of $S$, and where
    $Mor_{\mathcal{F}_S(G)}(P,Q)=Hom_G(P,Q)=\{c_g\in Hom(P,Q)|g\in G, gPg^{-1}\leq Q\}
    \cong N_G(P,Q)/C_G(P)$.
    
		Here $c_g$ denotes the homomorphism \textbf{conjugation} by $g$ $(x\mapsto gxg^{-1})$.
For each $P\leq S$, let $C'_G(P)$ be the maximal $p-$perfect subgroup of $C_G(P)$.Let $\mathcal{F}$ be a fusion system over a $p$-group $S$. For a discrete group $G$ and any finite set $\mathcal{H}$ of finite subgroups of $G$, let $\mathcal{T}_{\mathcal{H}}(G)$ denote the \textbf{transporter category of $G$}: the category with $Ob(\mathcal{T}_{\mathcal{H}}(G))=\mathcal{H}$, and where for each $P,Q\in\mathcal{H}$, $Mor_{\mathcal{T}_{\mathcal{H}}(G)}(P,Q)=N_G(P,Q)=\{g\in G|gPg^{-1}\leq Q\}$ the \textbf{transporter set}. The category $\mathcal{T}^c_{S}(G)$ is the full subcategory of $\mathcal{T}_{S}(G)$
  with objects the $\mathcal{F}_S(G)$-centric subgroups of $S$.
 A \textbf{linking system} $\mathcal{L}$
associated with $\mathcal{F}$ is a finite category together
with a pair of functors $\mathcal{T}_{Ob(\mathcal{L})}(S)\overset{\delta}{\longrightarrow}\mathcal{L}\overset{\pi }{\longrightarrow}\mathcal{F}$ such that the following conditions are satisfied: $Ob(\mathcal{L})$ is a set of subgroups of $S$ closed under $\mathcal{F}$-conjugacy and overgroups, and includes all subgroups which are $\mathcal{F}$-centric and $\mathcal{F}$-radical. Also $\delta$ is the
identity on objects, and $\pi$ is the inclusion on objects. and surjective on morphisms. For
each pair of objects $P,Q\in ob(\mathcal{L})$ such that $P$ is fully centralized in $\mathcal{F}$, $C_S(P)$ acts freely on
$Mor_{\mathcal{L}}(P,Q)$ via $\delta _{P,P}$ by right composition and $\pi _{P,Q}$ induces
a bijection
$Mor_{\mathcal{L}}(P,Q)/C_S(P)\overset{\simeq}{\longrightarrow}Hom_{\mathcal{F}}(P,Q).$
For each $P,Q\in ob(\mathcal{L})$ and each $g\in N_S(P,Q)$,
$\pi _{P,Q}$ sends $\delta _{P,Q}\in Mor_{\mathcal{L}}(P,Q)$ to $=c_g \in Hom_{\mathcal{F}}(P,Q)$. For all $\phi\in
Mor_{\mathcal{L}}(P,Q)$ and all $g\in P, \phi\circ\delta_{P,P}(g)=\delta
_{Q,Q}(\pi (\phi (g)))\circ \phi$.
A \textbf{centric linking system}
associated with a saturated fusion system $\mathcal{F}$ is a linking system with objects the set of $\mathcal{F}$-centric subgroups of $S$.
For every saturated fusion system there exists one associated centric linking system up to isomorphism \cite{Chermak}. 
 Let $\mathcal{L}^c_S(G)$ be the category whose objects are the $\mathcal{F}_S(G)-$centric subgroups of $S$, and where $Mor_{\mathcal{L}^c_S(G)}(P,Q)=N_G(P,Q)/C_G'(P).$
    Let $\pi:\mathcal{L}^c_S(G)\rightarrow\mathcal{F}_S(G)$ be the functor which is the inclusion on objects and sends the class of $g\in N_G(P,Q)$ to conjugation by $g$. For each $\mathcal{F}_S(G)-$centric subgroup $P\leq G$, let $\delta _P:P\rightarrow Aut_{\mathcal{L}^c_S(G)}(P)$ be the monomorphism induced by the inclusion $P\leq N_G(P)$.
A triple $(S,\mathcal{F},\mathcal{L})$ where $S$ is a finite
$p-$group, $\mathcal{F}$ is a saturated fusion system on $S$, and
$\mathcal{L}$ is an associated centric linking system with
$\mathcal{F}$, is called a $p-$local finite group. Its
\textbf{classifying space} is $|\mathcal{L}|\pcom$ where $(-)\pcom$
denotes the $p-$completion functor in the sense of Bousfield and
Kan and $|-|$ the classifying space functor. A space X is called \textbf{$p-$good} if the natural map $H_*(X;\mathbb{F}_p)\rightarrow H_*(X\pcom ;\mathbb{F}_p)$ is an isomorphism. Examples of spaces which are $p-$good are classifying spaces of finite groups. A finite group $G$ gives rise canonically to a $p-$local
finite group $(S,\mathcal{F}_S(G),\mathcal{L}^c_S(G))$ and
$BG\pcom\simeq |\mathcal{L}^c_S(G)|\pcom$ \cite{BK}. In particular, all
fusion systems coming from finite groups are saturated. In analogy with the Cartan-Eilenberg theorem the \textbf{ cohomology of a saturated fusion system} is defined to be the inverse limit over the orbit category. Recall \cite{BLO2} there exists an isomorphism of unstable algebras between the cohomology of the classifying space of a $p-$local finite group and the cohomology of the fusion system. This motivates us to define the \textbf{cohomology of the $p$-local finite group} to be this inverse limit.
Let $\mathcal{F}$ be a fusion system on the the finite $p-$group
$S$. $\mathcal{F}$ is called \textbf{Alperin} fusion system if
there are subgroups $P_1,P_2,\cdots P_r$ of $S$ and finite groups
$L_1,\cdots ,L_r$ such that for each $i$, $N_S(P_i)\in Syl_p(L_i)$,
$\mathcal{F}_{N_S(P_i)}(L_i)$ is contained in $\mathcal{F}$ and
$\mathcal{F}$ is generated by all the $\mathcal{F}_{N_S(P_i)}(L_i)$. Every saturated fusion system is Alperin
\cite[Section 4]{BLO1}. The groups $L_i=Aut_{\mathcal{L}}(P_i)$ 
will also be denoted $L_{P_i}$. 
A ring homomorphism $f:A\rightarrow B$ is called
$F-$\textbf{monomorphism} in the sense of Quillen \cite{Quillen1} if
every element in the kernel is nilpotent and
$F-$\textbf{epimorphism} in the sense of Quillen if every element in the cokernel is
nilpotent and $f$ is called $F-$\textbf{isomorphism} in the sense of Quillen if it is both
$F-$monomorphism in the sense of Quillen and $F-$epimorphism in the sense of Quillen. A monomorphism of unstable modules $i:N\rightarrow M$ is a \textbf{strong $F-$isomorphism} if and only if for all $x\in M, x\neq 0$, there exists an integer $c$ which depends on $x$ such that $Sq_0^c(x)\in N$, and $Sq_0^c(x)\neq 0$. 
 We recall that for $S \leq H\leq G$, $H$ is said \textbf{to control $p$-fusion in $G$}, if pairs of tuples of elements of $S$ are conjugate in $H$ if they are conjugate in $G$, or equivalently if for all $p$-subgroups $P,Q \leq S$, $N_H(P,Q)/C_H(P)$ equals $N_G(P,Q)/C_G(P)$ as homomorphisms from $P$ to $Q$, regardless of whether $G$ is finite or not. 

\subsection{Fusion systems of spaces}
Fusion systems and linking systems can be defined for spaces. We review what we need from \cite{BLO4}.
\begin{Definition}[\cite{BLO4}, Definition 1.6.]
Fix a space $X$, a finite $p$-group $S$, a map $f:BS\rightarrow X$.
\begin{enumerate}
\item Define $\mathcal{F}_{S,f}(X)$ to be the category whose objects are the subgroups of $S$, and whose morphisms are given by
\begin{eqnarray*}
Hom_{\mathcal{F}_{S,f}(X)}(P,Q)=\{\phi\in Inj(P,Q)|f|_{BP}\simeq f|_{BQ}\circ B\phi\}
\end{eqnarray*}
for each $P,Q\leq S$.
\item Define $\mathcal{F}'_{S,f}(X)\subseteq\mathcal{F}_{S,f}(X)$ to be the subcategory with the same objects as $\mathcal{F}_{S,f}(X)$, and where $Mor_{\mathcal{F}'_{S,f}(X)}(P,Q)$ (for $P,Q\leq S$) is the set of all composites of restrictions of morphisms in $\mathcal{F}_{S,f}(X)$ between $\mathcal{F}_{S,f}(X)-$centric subgroups.
\item Define $\mathcal{L}^c_{S,f}(X)$ to be the category whose objects are the $\mathcal{F}_{S,f}(X)-$centric subgroups of $S$, and whose morphisms are defined by
\begin{eqnarray*}
Mor_{\mathcal{L}^c_{S,f}(X)}(P,Q)=\{(\phi, [H])|\phi\in Inj(P,Q), H:BP\times I\rightarrow X,|\\H|_{BP\times 0}= f|_{BP}, |H|_{BP\times 0}= f|_{BQ}\circ B\phi\}
\end{eqnarray*}
The composite in $\mathcal{L}^c_{S,f}(X)$ of morphisms $P\rightarrow Q\rightarrow R$, where $H:BP\times I\rightarrow X$ and $K:BQ\times I\rightarrow X$ are homotopies as described above, are defined by setting
\begin{eqnarray*}
(\psi, [K])\circ (\phi,[H])=(\psi\circ\phi, [K\circ (B\phi \times Id)\cdot H]),
\end{eqnarray*}
where $\cdot $ denotes composition (juxtaposition) of homotopies. Let
\begin{eqnarray*}
\pi :\mathcal{L}^c_{S,f}(X)\rightarrow\mathcal{F}_{S,f}(X)
\end{eqnarray*}
be the forgetful functor: it is the inclusion on objects, and sends a morphism $(\phi, [H])$ to $\phi$. For each $\mathcal{F}_{S,f}(X)-$centric subgroup $P\leq S$, let 
\begin{eqnarray*}
\delta _P:P\rightarrow Aut_{\mathcal{L}^c_{S,f}(X)}(P)
\end{eqnarray*}
be the "distinguished homomorphism" which sends $g\in P$ to $(c_g, [f|_{BP}\circ H_g])$.

\end{enumerate}
\end{Definition}
Recall that when $S$ is a $p$-group, a map $f:BS\rightarrow X$ will be called \textbf{Sylow} if every map $BP\rightarrow X$, for a $p$-group $P$, factors through $f$ up to homotopy. A map $f:X\rightarrow Y$ between arbitrary spaces will becalled \textbf{centric} if the induced map $Map(X,X)_{Id}\stackrel{f\circ -}{\longrightarrow} Map(X,Y)_f$ is a homotopy equivalence.

We need a criteria to characterize spaces wich have cohomology isomorphic to the stable elements.
\begin{Theorem}[\cite{BLO2}, Theorem 7.5.]
A space $X$
which is $p$-complete 
 is the classifying space of some $p$-local finite group if and only if there is a finite $p$-group $S$ and a map $f:BS\rightarrow X$ such
that.
\begin{enumerate}
\item The fusion system $\mathcal{F}_{S,f}(X)$ is saturated.
\item There is a homotopy equivalence $X\simeq |\mathcal{L}^c_{S,f}(X)|\pcom$.
\item The restricted map $f|_{BP}$ is a centric map for each $\mathcal{F}_{S,f}(X)-$centric subgroup $P\leq S$.
\end{enumerate}
When these hold, $\mathcal{L}^c_{S,f}(X)$ is a centric linking system associated to $\mathcal{F}_{S,f}(X)$.

\end{Theorem}

\begin{Theorem}[\cite{BLO4}, Theorem 2.1.]
Let $X$ be a space, $S$ a finite $p$-group, $f:BS\rightarrow X$ a map. Assume that
\begin{enumerate}
\item $f$ is Sylow,
\item $f|_{BP}$ is a centric map for each $\mathcal{F}_{S,f}(X)-$centric subgroup $P\leq S$ and
\item every $\mathcal{F}'_{S,f}(X)-$centric subgroup of $S$ is also $\mathcal{F}_{S,f}(X)-$centric.
\end{enumerate}
Then the triple $(S,\mathcal{F}'_{S,f}(X),\mathcal{L}^c_{S,f}(X))$ is a $p$-local finite group.
\end{Theorem}
The obstruction theory in terms of characteristic classes for the realizabilty of the cohomology of a space was developed by Blanc \cite{Blanc}. This obstruction theory does not simplify in our situation.
\begin{Definition}
For a field $R$, let $coAlg_R$ denote the category of graded coalgebras over $R$.
\end{Definition}
\begin{Definition}
There is a sequence of characteristic 
cohomology
classes $\chi _n\in H^{n+2}(K_*;\Sigma ^nK_*)$
such that $\chi _n$ is defined whenever $\chi _1=\chi _2=\cdots =\chi _{n-1}=0$.
\end{Definition}
\begin{Theorem}[\cite{Blanc}, Proposition 6.5]
For $R=\mathbb{F}_p$, an unstable coalgebra $K_*\in \mathcal{C}A_R$ of finite projective dimension is realizable as the homology of some space $X\in\mathcal{T}_*$ if and only if $K_*$ has a coherently vanishing sequence of characteristic 
cohomology
classes $\chi _1=\chi _2=\cdots =\chi _{n-1}=0$.

\end{Theorem}
\subsection{Control of fusion in finite groups}
In finite groups fusion is controlled via $F-$isomorphism as shown by Benson, Grodal, Henke \cite{BGH}.
\begin{Theorem}[Theorem A, \cite{BGH}]
Let $\iota :H\leq G$ be an inclusion of finite groups of index prime to $p$, $p$ an odd prime, and consider the induced map on mod  $p$ group cohomology $\iota ^*:H^*(G;\mathbb{F}_p)\rightarrow H^*(H;\mathbb{F}_p)$. If for each $x\in H^*(H;\mathbb{F}_p)$, we have $x^{p^k}\in im(\iota ^*)$ for some $k\geq 0$, then $H$ controls $p$-fusion in $G$.
\end{Theorem}

\begin{Theorem}[Theorem B, \cite{BGH}]
Let $\mathcal{G}\leq\mathcal{F}$ be two saturated fusion systems on the same finite $p-$group $S$. Suppose that
$Hom_{\mathcal{G}}(A,B)=Hom_{\mathcal{F}}(A,B)$ for all $A,B\leq S$ with $A,B$ elementary abelian if $p$ is odd, and abelian of exponent at most $4$ if $p=2$. Then $\mathcal{G}=\mathcal{F}$.
\end{Theorem}
\subsection{Graph of groups}
A \textbf{finite directed graph} $\Gamma$
consists of two sets, the \textbf{vertices} $V$ and the
\textbf{directed edges} $E$, together with two functions $\iota,
\tau : E\rightarrow V$. For $e \in E,\iota (e)$ is called the
\textbf{initial vertex} of $e$ and $\tau (e)$ is the
\textbf{terminal vertex} of $e$. Multiple edges and loops are
allowed in this definition. The graph $\Gamma$ is \textbf{connected}
if the only equivalence relation on $V$ that contains all $(\iota
(e),\tau (e))$ is the relation with just one class. A graph $\Gamma$
may be viewed as a category, with objects the disjoint union of $V$
and $E$ and two non-identity morphisms with domain $e$ for each
$e\in E$, one morphism $e\rightarrow \iota (e)$ and one morphism
$e\rightarrow \tau (e)$. A \textbf{graph $\Gamma$ of groups} is a
connected graph $\Gamma$ together with groups $G_v, G_e$ for each
vertex and edge and injective group homomorphism
$f_{e,\iota}:G_e\rightarrow G_{\iota (e)}$ and $f_{e,\tau
(e)}:G_e\rightarrow G_{\tau (e)}$ for each edge $e$. 
Let $(S,\mathcal{F},\mathcal{L})$ be a $p$-local finite group, $P_1,\dots,P_n$ a collection of $\mathcal{F}$-centric subgroups of $S$ which might contain multiples, with a generating collection $K=\{K_1,\dots,K_n\}$ where $K_i<\Aut_\mathcal{L}(P_i)$. We define a \textbf{graph of groups $\Gamma_K$ associated with $K$}. 
%
%
Let $F_{i,j}$ be the subgroup of $\Aut_\mathcal{F}(\langle P_i,P_j\rangle)$ consisting of automorphisms that restrict in $\mathcal{F}$ to both an automorphism of $P_i$ contained in $\pi(K_i)$ and an automorphism of $P_j$ contained in $\pi(K_j)$. Define $K_{i,j}$ as the preimage of $F_{i,j}$ in $\Aut_\mathcal{L}(\langle P_i,P_j\rangle)$. Note that $K_{i,j}=K_{j,i}$. 
By \cite[Prop 2.11.]{Lib06} we can define an injective restriction homomorphism $k_{i,j}: K_{i,j}\rightarrow K_j$.
 A graph of groups $\Gamma_K$ will be called a graph \textbf{associated} with $K$, if there exist Sylow $p$-subgroups $H_i\leq K_i$ and $H_{i,j}<K_{i,j}$ such that $\Gamma_K$ has $n$ vertices, with corresponding groups $K_i$ and such that there is a directed edge from $K_i$ to $K_j$ in $\Gamma_K$ if and only if $H_j<k_{i,j}(H_{i,j})$, and to this edge correspond the injective group homomorphisms $k_{j,i}:  K_{i,j}\rightarrow K_i$ and $k_{i,j}:  K_{i,j}\rightarrow K_j$.
In our applications we will need to choose a special subgraph of $\Gamma_K$.
An \textbf{aborescence} is a directed graph in which, for a vertex $u$ called the \textbf{root} and any other vertex $v$, there is exactly one directed path from $u$ to $v$.
 Let $\Gamma_K$ be a graph associated with $K$. An arborescence $T$ with root $t$ of $\Gamma_K$ is called a \textbf{generating tree of $\Gamma_K$}, if $P_t\vartriangleleft S$ and $H_t=S<K_t$.
 \subsection{Group models for fusion systems and their representations}
 Let $\mathcal{F}$ be a saturated fusion system on a finite $p-$group $S$.
 A discrete group $\mathcal{G}$ is a \textbf{group model for
$\mathcal{F}$} if $S$ is a Sylow $p-$subgroup of $\mathcal{G}$ and
$\mathcal{F}_S(\mathcal{G})=\mathcal{F}$ in the sense that all morphisms of $\mathcal{F}$ are induced by conjugation of elements in $\mathcal{G}$. A saturated fusion system for which exists no finite group model is called \textbf{exotic}. A \textbf{signalizer functor} on a group model $\mathcal{G}$ is an assignment $P\mapsto\theta (P)$ for every $\mathcal{F}-$centric subgroup $P\leq S$ such that $\theta (P)$ is a complement of $Z(P)$ in $C_{\mathcal{G}}(P)$ and such that if $gPg^{-1}\leq Q$ for $g\in\mathcal{G}$ then $\theta (Q)\leq g\theta (P)g^{-1}$. A signalizer functor if it exists gives rise to a centric linking system \cite{LS}. All group models for fusion systems known so far \cite{Ian+Radu}, \cite{Robinson1}, \cite{LS}, \cite{anewfam} have signalizer functors \cite{LS}, \cite{anewfam}. Aschbacher and Chermak \cite{AschbacherChermak} introduced the notion of representation of a $p$-local finite group if there exists a group model for $\mathcal{F}$ and a signalizer functor on this group which induces $\mathcal{L}$. All so far known group models \cite{Ian+Radu}, \cite{Robinson1}, \cite{LS}, \cite{anewfam} have such a representation as shown by ourselves in two independent collaborations one together with Libman \cite{LS} and in \cite{anewfam}.\\
In \cite{gmffs} we prove that for every group model we have a map in $\mathbb{F}_p-$cohomology into the stable elements $H^*(B\mathcal{G})\overset{q}{\rightarrow} H^*(\mathcal{F})$. 
A discrete group $\mathcal{G}$ is a \textbf{group model for
$(S,\mathcal{F},\mathcal{L})$} if it is represented in the sense of Aschbacher and Chermak and has the $\mathbb{F}_p-$cohomology of $(S,\mathcal{F},\mathcal{L})$.  For constrained fusion systems there exists finite group models \cite{BCGLO1} which are therefore group models for $(S,\mathcal{F},\mathcal{L})$. The question whether for every $p$-local finite group $(S,\mathcal{F},\mathcal{L})$ there exists a group model is open.
We review the existing group model constructions \cite{Ian+Radu}, \cite{Robinson1}, \cite{LS}, \cite{anewfam} for fusion and linking systems.

The first group model constructed by Leary and Stancu consists of an iterated HNN-construction.
\begin{Theorem}[\cite{Ian+Radu}, Theorem 2]
Let $\mathcal{F}$ be a fusion system on $S$ generated by $\Phi=\{\phi_1, \cdots, \phi_r\}$. Let $T$ be a free group with free generators $t_1, \ldots, t_r$, and define $\pi_{LS}$ as the quotient of the free product $S*T$ by the relations $t_i^{-1}ut_i=\phi_i(u)$ for all $i$ and for all $u\in P_i$. Then $\pi_{LS}$ is a group model for $\mathcal{F}$. 
\end{Theorem}

Group models for fusion systems involving only amalgams
(for example Robinson \cite{Robinson1}, Libman and ourselves  \cite{LS}
)
are a special case of the following construction due to ourselves and Leip \cite{anewfam}.
\begin{Theorem}[\cite{anewfam}, Theorem 4.1.]
\label{groupmodelkt}
 Let $K=\{K_1,\dots,K_n\}$ be a generating collection, $\Gamma_K$ the associated graph for some choice of $H$'s and $T$ a generating tree in $\Gamma_K$. Then the amalgam $\pi_{K,T}$ over the graph of groups $T$ is a group model for $S$ and a group model for $\mathcal{L}$.
\end{Theorem}
The finite group $G=S\wr
\Sigma _{e(X)}$ constructed by Park \cite{Sejong} is never a group model since the group $S$ is never a Sylow $p-$subgroup of $G$, where $e(X)$ is a certain characteristic biset associated with $\mathcal{F}$.\\

A group model is \textbf{minimal} if the map in the following theorem is an isomorphism in $\mathbb{F}_p-$cohomology.
\begin{Theorem}[\cite{gmffs}]
Let $\mathcal{F}$ be a saturated fusion system over the finite
$p-$group $S$ and $\mathcal{G}$ a group model for $\mathcal{F}$.
Then there exist a natural map of unstable algebras
$H^*(B\mathcal{G})\overset{q}{\rightarrow} H^*(\mathcal{F})$ making
$H^*(\mathcal{F})$ a module over $H^*(B\mathcal{G})$.
\end{Theorem}
The cohomology ring $H^*(B\mathcal{G};\mathbb{F}_p)$ is isomorphic to the stable elements up to $F-$isomorphism \cite{gmffs}.
\begin{Theorem} Let $\mathcal{F}$ be a saturated fusion system over the finite $p-$group 
$S$ and $\mathcal{G}$ one of the above group models for $\mathcal{F}$. Then there exist natural maps of algebras over the Steenrod
algebra $q: H^* (B\mathcal{G}) \rightarrow  H^* (\mathcal{F})$ and $r^*: H^*(\mathcal{F}) \rightarrow H^*(B\mathcal{G})$ such that we obtain a split
short exact sequence of unstable modules 
$0\rightarrow W \rightarrow H^* (B\mathcal{G}) \overset{\leftarrow}{\rightarrow} H^* (\mathcal{F}) \rightarrow 0,$
with
$W \cong Ker(Res^{\mathcal{G}}_S )$ as unstable module.
\end{Theorem}
The fact that $W$ is not always trivial illustrates that the Martino-Priddy Conjecture does not hold in the general case. In fact we have that.
The question of existence of minimal group models for saturated fusion systems is related to the existence of the centric linking systems in the following way. 
 Recall \cite{BLO4} that an associated centric linking system exists if there exists a space which has the homology of the linking system. This motivates our continued search for minimal group models.\\ 

\subsection{Fusion systems over discrete $p$-toral groups and over profinite groups}
A group $G$ is \textbf{locally finite} if every finitely generated subgroup of $G$ is finite, and is a \textbf{locally finite $p$-group} if every finitely generated subgroup is a finite $p$-group, and is \textbf{artinian} if every nonempty set of subgroups, partially ordered by inclusion, has a minimal element. 
Let $\mathbb{Z}/p^{\infty}$ be the union of the cyclic groups $Z/p^n$ under the canonical inclusions. A \textbf{discrete $p$-toral group} is a group $P$ which contains a normal subgroup $T_p\leq P$, isomorphic to a finite product of copies of $\mathbb{Z}/p^{\infty}$, and such that the quotient $P/T_p$ is a finite $p$-group. For such a group, we say that $T_p$ is the \textbf{connected component} or \textbf{maximal torus} of $P$ and $P/T_p$. The \textbf{rank} of a discrete $p$-toral group $P$ is the rank of the maximal torus $T_p$. That is, if $T_p\cong (\mathbb{Z}/p^{\infty})^k$, then we say that $rk(P)=k$.
A  \textbf{fusion system} $\mathcal{F}$ over a discrete $p$-toral group $S$ is a category whose objects are all the subgroups of $S$, and which satisfies the following two properties for all $P,Q\leq S$: The set $Hom_{\mathcal{F}}(P,Q)$ contains injective group homomorphisms and amongst them all morphisms induced by conjugation of elements in $S$ and each element is the composite of an isomorphism in $\mathcal{F}$ followed by an inclusion. Two subgroups $P,Q\leq S$ will be called $\mathcal{F}-$conjugate if they are isomorphic in $\mathcal{F}$.
Let $\mathcal{F}$ be a fusion system over a discrete $p$-toral group $S$. 
 A subgroup $P\leq S$ is \textbf{$\mathcal{F}$-fully centralized} resp. \textbf{fully $\mathcal{F}$-normalized} if $|C_S(P)|\geq |C_S(P')|$ resp. $|N_S(P)|\geq |N_S(P')|$ for all $P'\leq S$ which is $\mathcal{F}$-conjugate to $P$.
 Let $\mathcal{F}$ be a saturated fusion system over a discrete $p$-toral group $S$. A \textbf{fusion-controlling set} is a set $\mathcal{P}=\{P_1,\cdots, P_r\}$ of representatives of the $\mathcal{F}-$conjugacy classes of $\mathcal{F}-$centric for all
 $\mathcal{F}-$radical subgroups such that, for each $j$, $P_j$ is fully $\mathcal{F}-$normalized.
In fact, it is sufficient to restrict oneself to representatives of the
essential subgroups.
\\
 A fusion system $\mathcal{F}$ is called \textbf{saturated} if for all $P\leq S$ which is fully normalized in $\mathcal{F}$, $P$ is fully centralized in $\mathcal{F}$ and $Aut_S(P)\in Syl_p(Aut_{\mathcal{F}}(P))$ and moreover if $P\leq S$ and $\phi\in Hom_{\mathcal{F}}(P,S)$ are such that $\phi (P)$ is fully centralized, and if we set
$N_{\phi}=\{g\in N_S(P)|\phi c_g \phi^{-1}\in Aut_S(\phi (P))\}$,
then there is $\overline{\phi}\in Hom_{\mathcal{F}}(N_{\phi},S)$ such that $\overline{\phi}|_P=\phi$. A subgroup $P\leq S$ will be called $\mathcal{F}-$centric if $C_S(P')\leq P'$ for all $P'$ which are $\mathcal{F}-$conjugate to $P$. Denote $\mathcal{F}^c$ the full subcategory of $\mathcal{F}$ with objects the $\mathcal{F}-$centric subgroups of $S$.  A saturated fusion system $\mathcal{F}$ over a discrete $p$-toral group $S$ is \textbf{constrained} if it contains an object which is $\mathcal{F}$-centric and $\mathcal{F}$-normal.
A \textbf{centric linking system associated with $\mathcal{F}$} is a category $\mathcal{L}$ whose objects are the
$\mathcal{F}$-centric subgroups of $S$, together with a functor
$\pi :\mathcal{L}\longrightarrow\mathcal{F}^c$,
and \textbf{"distinguished" monomorphisms $\delta _P:P\rightarrow Aut_{\mathcal{L}}(P)$} for each $\mathcal{F}$-centric subgroup $P\leq S$ such that the following conditions are satisfied: $\pi$ is the identity on objects and surjective on morphisms. More precisely, for each pair of objects $P,Q\in\mathcal{L},Z(P)$ acts freely on $Mor_{\mathcal{L}}(P,Q)$ by composition (upon identifying $Z(P)$ with $\delta _P (Z(P))\leq Aut_{\mathcal{L}}(P)$), and $\pi$ induces a bijection
$Mor_{\mathcal{L}}(P,Q)/Z(P)\overset{\simeq}{\longrightarrow}Hom_{\mathcal{F}}(P,Q).$
For each $\mathcal{F}$-centric subgroup $P\leq S$ and each $x\in P$, $\pi (\delta _P(x))=c_x \in Aut_{\mathcal{F}}(P)$.
For each $f\in Mor_{\mathcal{L}}(P,Q)$ and each $x\in P, f\circ\delta_P(x)=\delta _Q(\pi f(x))\circ f$.
A \textbf{$p$-local compact group} is a triple $(S,\mathcal{F},\mathcal{L})$, where $S$ is a discrete $p$-toral group, $\mathcal{F}$ is a saturated fusion system over $S$, and $\mathcal{L}$ is a centric linking system associated to $\mathcal{F}$. The \textbf{classifying space} of $\mathcal{G}$ is the $p-$completed nerve $B\mathcal{G}=|\mathcal{L}|\pcom$. Levi and Libman \cite{LL} showed that for every saturated fusion system over a discrete $p$-toral group there exists a unique classifying space up to isomorphism. In his thesis [\cite{Gonzalez}, Theorem 2.2.3] Gonzalez proved that there exists a stable elements for fusion system over discrete $p$-toral groups for $\mathbb{F}_p-$cohomology under the assumption that the $p$-local compact group can be approximated by $p$-local finite groups. He later extended his result and showed this is true in the general case and for coefficients in a trivial $\mathbb{Z}_{(p)}-$module, [\cite{GonzalezStableElements}, Theorem 4.4]. We expect that this result 
will prove
to be true for twisted coefficients as well.\\
A \textbf{profinite group} is a group $G$ that is an inverse limit of finite groups and made into a topological space using the profinite topology. They have Sylow $p$-subgroups and Stancu and Symonds construct fusion systems over profinite groups first using pro-fusion systems over pro-$p$-groups \cite{StancuSymonds}. 

\subsection{The Robinson model for fusion systems over discrete $p$-toral groups}
For every fusion system $\mathcal{F}$ over a discrete $p$-toral group $S$ there exists an infinite group $\mathcal{G}$ such that $\mathcal{F}_S(\mathcal{G})=\mathcal{F}$.
\begin{Proposition}[\cite{Gonzalez}, Prop. 2.1.3]
Let $\mathcal{F}$ be a constrained saturated fusion system over a discrete $p$-toral group $S$. There exists a unique up to isomorphism $p'$-reduced  $p$-constrained artinian locally finite group $G$ such that $G$ is a group model for $\mathcal{F}$ and if $\mathcal{L}$ is the centric linking system associated with $\mathcal{F}$ then $G\cong Aut_{\mathcal{L}}(P)$ for any $\mathcal{F}$-centric subgroup $P$ which is normal in $\mathcal{F}$, $\mathcal{L}\cong \mathcal{L}_S(G)$.
\end{Proposition}

The groups of Robinson type are iterated amalgams of automorphism groups in the linking system over the $S-$normalizers of the respective $\mathcal{F}$-centric subgroups of $S$.  Note that these automorphism groups exist and are known regardless of whether $\mathcal{L}$ exists or not.
In his thesis, Gonzalez \cite{Gonzalez} proved an analogue of the Robinson model \cite{Robinson1} for discrete $p$-toral groups which we will generalize.
\begin{Theorem}[\cite{Gonzalez}, Theorem 2.2.3.]
Let $\mathcal{F}$ be a saturated fusion system over a discrete
$p-$toral group $S$ and let $\mathcal{P}=\{P_1=S, P_2,\cdots, P_r\}$ be a fusion-controlling set for $\mathcal{F}$. Furthermore, for each $P_j\in\mathcal{P}$, let $L_j$ be the group obtained above.
Then
\begin{eqnarray*}
\mathcal{G}=L_1\underset{N_S(P_2)}{*}L_2\underset{N_S(P_3)}{*}...\underset{N_S(P_n)}*L_n
 \end{eqnarray*}
 is a group model for $\mathcal{F}$.
\end{Theorem}

\subsection{The Kan-Thurston Theorem}
Kan and Thurston \cite[Theorem 1.1]{KanThurston} proved the following result which we need at a later stage.
\begin{Theorem}
Let $X$ be a path connected space with base point. There
exists
a
fibration
in the sense of Serre
$tX: TX\rightarrow X$
which
is
natural
with
respect
to
$X$
and
has
the
following
properties.
\begin{enumerate}
\item
The
map
$tX$
induces
an
isomorphism
on
singular
homology
and
cohomology
with
local
coefficients
$H_*(TX;A)\cong
H_*(X;A)$ and $
H^*(TX;A)\cong
H^*(X;A)$
for
every
local
coefficient
system
$A$
on
$X,$
and
\item
$\pi _iTX$ is trivial for $i\neq 1$ and $\pi _1tX$ is onto.
\end{enumerate}
Furthermore
the
homotopy
type
of
$X$
is
completely
determined
by
the
pair
of
groups
$(Gx,
Px)$
where
$Gx
=\pi _1TX$
and
$Px
=
ker\pi _1tX.$
\end{Theorem}
\subsection{Unstable modules over the Steenrod Algebra and Lannes $T-$Functor}
The category $\mathcal{U}$ of unstable modules over the Steenrod algebra admits a decreasing filtration by full subcategories $\mathcal{N}il_{k}$, $k\geq 0$. These categories are defined as follows. The catgory $\mathcal{N}il_{k}$ is the smallest full subcategory, stable under direct limits and containing $\Sigma ^kM$ for every unstable module $M$. Every unstable module admits a decreasing filtration which will be called \textbf{nilpotent filtration} of the module. If $M$ is an unstable its biggest unstable submodule $s-$nilpotent will be denoted $M_s$. The quotient $M_s/M_{s+1}$ is of the form $\Sigma ^sR_s(M)$ where $R_s(M)$ is a reduced unstable module meaning it does not contain any nontrivial suspension. 
Two unstable reduced modules $M$ and $N$ are called \textbf{strongly $F-$isomorphic} if they have the same $\mathcal{N}il-$localisation and will be denoted $M\cong _{F}N$. The $n-$th \textbf{Brown-Gitler module} $J(n)$ is the representing functor for $H_n$, $H_n(M)\cong Hom_{\mathcal{U}}(M,J(n))$. As usual $\mathcal{K}$ denotes the category of unstable algebras over the Steenrod algebra.\\[0.3cm]
Lannes and Schwartz characterize nilpotent modules in terms of the injective objects $H^*(BV;\mathbb{F}_p)$. 
\begin{Theorem}
An unstable module $N$ is nilpotent if and only if
$ Hom_{\mathcal{U}_n}(N, H^*(BV;\mathbb{F}_p))=0$ for all elementary abelian $p$-groups $V$.
\end{Theorem}

For an elementary abelian $p$-group $V$, Lannes considered the functor $T_V:\mathcal{U}\rightarrow\mathcal{U}$ characterized by the property $Hom_{\mathcal{U}}(T_VM,N)\cong Hom_{\mathcal{U}}(M, H^*BV\otimes N)$ which is exact and commutes with tensor products. The functor $T_V$ lifts to a functor from $\mathcal{K}$ to itself and continues to be adjoint to $H^*(BV;\mathbb{F}_p)$.
An $\mathcal{A}-$module $M$ is said to be \textbf{locally finite} if it is direct limit of finite modules or equivalently iff, for any $x$ in $M$, the span of $x$ over $\mathcal{A}_p$ is finite as a set or as an $\mathbb{F}_p-$vector space. 
\begin{Theorem}[\cite{Schwartz}, Theorem 6.2.1] Let $M$ be an unstable $\mathcal{A}_p
$-module. It is equivalent that:
\begin{enumerate}
\item $M$ is locally finite;
\item $Hom_{\mathcal{U}}(M,\overline{H}^*V\otimes J(m))$ is trivial for any $V$ and $m\geq 0$;
\item $T_VM$ is isomorphic to $M$ for any $V$;
\item $TM$ is isomorphic to $M$;
\item $\overline{T}_VM$ is trivial for any $V$;
\item $\overline{T}M$ is trivial.
\end{enumerate}
\end{Theorem}
\begin{Proposition}
Let $M$ be an unstable module and $k$ a fixed integer, $k\geq 1$. Assume we have a short exact sequence: $0\rightarrow K\rightarrow M\rightarrow N\rightarrow 0$, with $K\in\mathcal{N}il_k$. If $s<k$ we have $R_s(M)\cong R_s(N)$.
\end{Proposition}
\begin{Theorem}
Let $G$ be a finite group and $V$ be an elementary abelian $p$-group. Choose a representative for each element $\rho$ in 
$Rep(V,G)$. Then the resulting map \begin{eqnarray*}
l^V_G:T_VH^*BG\rightarrow \underset{\rho\in Rep(V,G)}{\prod{H^*BC_G(\rho)}}
\end{eqnarray*}
 with components $ad(c^*_{\rho})$ is an isomorphism in $\mathcal{K}$.
\end{Theorem}
\begin{Definition}
A topological group $G$ is a \textbf{Lannes group} at the prime $p$ if the map $l^V_G$ is an isomorphism for all elementary abelian $p$-groups $V$ (with $H^*(BG;\mathbb{F}_p)$ and $H^*BC_G(\rho)$ replaced by the continuous mod $p$-cohomology $H^*_{cts}(G;\mathbb{F}_p)$ and is assumed to be noetherian if $G$ is profinite).
\end{Definition}
\subsection{Universally stable elements}
Let $G$ be a finite group with a subgroup $P$. An element $x\in H^*(P;R)$ is called \textbf{stable} if for all subgroups $H,Q\leq P$ and elements $g\in G$ such that $gHg^{-1}\leq Q $ we have $Res \circ c_g^*(x)=Res(x)$. The Cartan-Eilenberg Theorem \cite[Theorem XII 10.1]{CE} states that the cohomology ring of a finite group consists of the ring of stable elements of the cohomology ring $H^*(P;R)$ of a $p-$Sylow subgroup $P$ for a solid ring of coefficients $R$. If an element in the cohomology ring $x$ is not stable it is called \textbf{unstable}.\\
Leary, Schuster, Yagita \cite{LearySchusterYagita} generalize the stable elements theorem of Cartan and Eilenberg \cite{CE}. For a finite $p$-group $P$, let $\mathcal{C}_u$ be the category whose objects are the subgroups of $P$, with morphisms all injective group homomorphisms. Let $\mathcal{C}$ be any subcategory of $\mathcal{C}_u$ such that $P$ is an object of $\mathcal{C}$, and such that for any object $Q$ of $\mathcal{C}$, the inclusion of $Q$ in $P$ is a morphism in $\mathcal{C}$. Let $H^*(P)$ stand for mod-$p$ group cohomology, which may be viewed as a contravariant functor from $\mathcal{C}_u$ to $\mathbb{F}_p-$algebras. We shall study the limit $I(P,\mathcal{C})$ of this functor: $I(P,\mathcal{C};R):=\underset{Q\in \mathcal{C}}{lim }H^*(Q;R)\subset H^*(P;R)$ for $R$ solid. 
Note that the cohomology ring $H^*(P;R)$ can be filtered by all the different limits of this functor.
\begin{Theorem}[\cite{LearySchusterYagita}, Theorem 1]
\label{LSY}
There exists a discrete group $\Gamma$ containing $P$ as a subgroup such that:
\begin{enumerate}
\item $Im(Res^{\Gamma}_P)$ is equal to $I(P,Q)$;
\item $(Ker(Res^{\Gamma}_P))^2$ is trivial;
\item $Res^{\Gamma}_P$ induces an isomorphism from $H^*(\Gamma)/\sqrt{0}$ to $I(P,\mathcal{C})/\sqrt{0}$;
\item $\Gamma$ is virtually free.
\end{enumerate}
\end{Theorem}
The conditions on $\mathcal{C}$ can be relaxed as follows. Consider arbitrary finite categories with objects finite groups and morphisms injective group homomorphisms.
Define  $I(\mathcal{C})$ to be the limit over this category and for any group $\Gamma $, define $\mathcal{D}(\Gamma)$ to be the category of finite subgroups of $\Gamma $, with morphisms inclusions and conjugation by elements of $\Gamma $.
\begin{Theorem}[\cite{LearySchusterYagita}, Theorem 1']
Let $\mathcal{C}$ be a connected finite category of finite groups and injective homomorphisms. Then there exists a discrete group $\Gamma$ and a natural transformation from $\mathcal{C}$ to $\mathcal{D}(\Gamma )$ such that $\Gamma $ and the induced map from $H^*(\Gamma)$ to $I(\mathcal{C})$ satisfy properties $1.)$, $2.)$, $3.)$, $4.)$ of the previous Theorem.
\end{Theorem}
\subsection{Mislin's Theorem for fusion systems}
There is control of fusion via cohomology in the analogue of Mislin's theorem for fusion systems.
\begin{Theorem}[\cite{Park3}, Theorem 1] Let $k$ be a field of characteristic
$p$. Let $\mathcal{F}$ be a saturated fusion system on a finite $p$-group $S$ and let $\mathcal{E}$ be a
saturated subsystem of $\mathcal{F}$ on the same $p$-group $S$. Suppose that $H^*
(\mathcal{F}; k) \cong H^*
(\mathcal{E}; k)$.
Then $\mathcal{F} = \mathcal{E}$.
\end{Theorem}
\subsection{Fusion and complex oriented cohomology theories}
We review what we need from \cite{Schuster}. A multiplicative generalised cohomology theory $E$ is called \textbf{complex oriented}, if
there is a class $x = x^E
 \in E^2
(CP ^{
\infty })$, called a \textbf{complex orientation}, which pulls
back to a generator of the free rank one module $E^2(CP^1)$ over the respective cohomology theory $E(-)$under the canonical inclusion
$CP^
1 \subset CP^{\infty }
$. Examples are singular cohomology $HR$ for a commutative ring $R$, complex $K-$theory, complex cobordism $MU$ with coefficients $MU_*$, the Brown-Peterson spectrum $BP$, the John Wilson spectra $E(n)$.
\subsection{Group models for fusion systems}
We quote a general construction.
\begin{Theorem}
\label{groupmodelkt}
 Let $K=\{K_1,\dots,K_n\}$ be a generating collection, $\Gamma_K$ the associated graph for some choice of $H$'s and $T$ a generating tree in $\Gamma_K$. Then the amalgam $\pi_{K,T}$ over the graph of groups $T$ is a group with the following properties:
\begin{enumerate}
 \item $S$ is a Sylow $p$-subgroup of $\pi_{K,T}$.
 \item $\mathcal{F}=\mathcal{F}_S(\pi_{K,T})$.
\item $H^*(|\mathcal{L}|,\F_p)$ is a retract of $H^*(B\pi_{K,T})$ in the category of unstable algebras. It is equal to the image of $H^*(\pi_{K,T},\F_p)\rightarrow H^*(S,\F_p)$, the product of any two elements
 in the kernel is zero.
 \item $B\pi_{K,T}$ is $p$-good.
 \item $H^*(B\pi_{K,T})$ is finitely generated.
 \item $|\mathcal{L}|\pcom$ is a stable retract of $(B\pi_{K,T} )\pcom $.
 \item $H^*(B\pi_{K,T})$ is $F-$isomorphic to the stable elements in
 the sense of Quillen.
\end{enumerate}
\end{Theorem}
\section{An illustration of the $F-$ isomorphism theorem of Quillen}
 It is a consequence of Quillen's theorem that for reasonable discrete groups the $\mathbb{F}_p-$cohomology is determined up to $F$-isomorphism by the stable elements. In more detail, a group $G$ is called a \textbf{Quillen group} if the map $q_G:H^*(BG;\mathbb{F}_p)\rightarrow lim_{\mathcal{A}_p(G)^{op}}H^*(BE;\mathbb{F}_p)$ induced by the restriction homomorphisms, where $\mathcal{A}_p(G)^{op}$ is the opposite Quillen category whose objects are all the elementary abelian $p-$subgroups of $G$, is an $F-$isomorphism. We denote the ring not unstable algebra $H^*(Q)=lim_{\mathcal{A}_p(G)^{op}}H^*(BE;\mathbb{F}_p)$. All group models known so far fulfill this condition which we illustrate here in detail. At a later stage we will show how it is possible to modify these group models such that the map $q_G$ is not an $F$-isomorphism anymore. The remaining situation is, what can be said about the cohomology of a space which has the right fusion system.\\
The first theorem of this section is a special case of the theorem we prove at the end of this section. The techniques of the proof are quite different and, as we believe, interesting in themselves so we include both. This is why we would like to include this special case separately. Results hold for coherent algebras \cite{Higman}.
\begin{Theorem} Let $S$ be a finite $p-$group and $\mathcal{F}$ a fusion system over $S$. Let $\mathcal{G}$ be a group model for $\mathcal{F}$ with finite virtual dimension and
 $H^*(B\mathcal{G};\mathbb{F}_p)$ Noetherian. Then $H^*(B\mathcal{G};\mathbb{F}_p)$ is $F-$monomorphic in the sense of Quillen to $H^*(\mathcal{F})$.\end{Theorem}

\underline{Proof:} From the theorem above we have a map $f:H^*(B\mathcal{G})\rightarrow H^*(\mathcal{F})$. We want to show that all elements in the kernel of $f$ are nilpotent. We will show that the nilpotent elements in $H^*(B\mathcal{G};\mathbb{F}_p)$ are precisely the ones who lie in the kernel of all the restriction maps to all the elementary abelian subgroups $V\leq \mathcal{G}$.
Let $A$ be a commutative graded connected algebra over the field $\mathbb{F}$ and $I\subset A$ an ideal. The set $Ass(I)$ of associated prime ideals of $I$ is finite and $\sqrt{I}=\bigcap \left\{p| p \in Ass(I) \right\} $.
Since ideals in a Noetherian algebra are finitely generated an ideal in which all elements are nilpotent is a nilpotent ideal. We will show that the kernel of the restriction map is contained in the nilradical $n$ of $H^*(B\mathcal{G};\mathbb{F}_p)$.
We will show that an element $x\in ker(res^*)$ belongs to every minimal prime ideal. 
Let $\mathcal{P}^*$ denote the Steenrod algebra of the Galois field $\mathbb{F}_p$ and recall that for any reasonable topological space $X$ that $H^*(X;\mathbb{F}_p)$ is an unstable algebra over $\mathcal{P}^*$. For any unstable algebra $H^*$ over $\mathcal{P}^*$ one says that an ideal $I\subset H^*$ is $\mathcal{P^*}-$invariant if $\theta (I)\subseteq I$ for all $\theta\in\mathcal{P}^*$. There is a $\mathcal{P}^*-$version of the Lasker-Noether-Theorem, \cite{NS}. Moreover we have that if $H^*$ is a Noetherian unstable algebra over the Steenrod algebra $\mathcal{P}^*$ and $I\subseteq H^*$ a $\mathcal{P}^*-$invariant ideal, then the associated prime ideals of $I$ are all $\mathcal{P}^*-$invariant.
This implies that the minimal primes of a Noetherian unstable $\mathcal{P}^*-$algebra are $\mathcal{P}^*-$invariant.
Let $H^*$ be a Noetherian unstable algebra which in our applications is going to be the coordinate ring. Denote by $min(H^*)$ the set of minimal prime ideals of $H^*$. For $p\in min (H^*)$ the quotient algebra  $H^*/p$ is therefore an unstable Noetherian integral domain over $\mathcal{P}^*$. By [1] we can therefore find an embedding $H^*/p\rightarrow \mathbb{F}_p[V_p]$ which is a finite ring extension, where $V_p$ is a finite dimensional vector space over $\mathbb{F}_p$ and $\mathbb{F}_p[V_p]$ is the polynomial algebra on $V_p$ regarded as an unstable $\mathcal{P}^*-$algebra in the canonical way. Moreover, if $\phi :H^*\rightarrow \mathbb{F}_p[V]$ is a map of unstable algebras then the kernel of $\phi$ is an invariant prime ideal, and if $\mathbb{F}_p[V]$ is finitely generated as an $H^*-$module even a minimal prime of $H^*$.
For a Noetherian cohomology algebra $H^*(BG;\mathbb{F}_p)$ this means the nil-radical can be described in the following way.
\begin{Theorem}
Let $H^*$ be a Noetherian unstable algebra over the Steenrod algebra. If $p\subset H^*$ is a minimal prime ideal then there is a homomorphism $\phi :H^*\rightarrow \mathbb{F}_p[V_{\phi}]$ into a finitely generated module over $H^*$. Therefore $h\in H^*$ is nilpotent if and only if $h$ belongs to the kernel of every homomorphism $\phi :H^*\rightarrow \mathbb{F}_p[V_{\phi}]$, where $V_{\phi}$ is a finite dimensional vector space over $\mathbb{F}_p$.
\end{Theorem} 
At this point we want to note there is a difference between the cases $p=2$ and $p\neq 2$. For $p=2$ one has $\mathbb{F}_2[V]\cong H^*(BV;\mathbb{F}_2)$, where as for $p\neq 2$ one only has $\mathbb{F}_p[V]\cong H^*(BV;\mathbb{F}_p)/\sqrt{0}$.  \\
We specialize the algebra $H^*$ in the previous theorem to  $H^*(BG;\mathbb{F}_p)$ where $G$ is some kind of reasonable discrete group and describe homomorphisms $\phi :H^*(BG;\mathbb{F}_p)\rightarrow \mathbb{F}_p[V_{\phi}]$ in terms of $p-$elementary abelian subgroups $V_{\phi}\leq G$ and the homomorphism induced by the inclusion. This leads to a number of technical problems that were nicely dealt with in \cite{DLS} and \cite{Zarati}. First of all we need to review some results of \cite{Quillen1} and \cite{Quillen2}. For this reason we are forced at first to introduce the full Steenrod algebra $\mathcal{A}^*$ of the Galois field $\mathbb{F}_p$, namely the algebra generated by not just the Steenrod reduced powers, but in addition to them also the Bockstein. From \cite{Quillen1} and \cite{Quillen2} we need the following fact (see \cite{DLS} Observation 2 for a short proof).
\begin{Theorem}[Quillen] If $H^*(BG;\mathbb{F}_p)$ is Noetherian and $\phi :H^*(BG;\mathbb{F}_p)\rightarrow H^*(BV;\mathbb{F}_p)$ is a homomorphism of algebras over the Steenrod algebra $\mathcal{A}^*$ making $H^*(BV;\mathbb{F}_p)$ into a finitely generated $H^*(BG;\mathbb{F}_p)-$module, then there is an inclusion $V\leq G$ such that the compositions
\begin{eqnarray*}
H^*(BG;\mathbb{F}_p)\underset{res^*}{\overset{\phi}{\overset{\longrightarrow}{\longrightarrow}}}H^*(BV;\mathbb{F}_p)\overset{\pi}{\rightarrow }H^*(BV;\mathbb{F}_p)/\sqrt{0}
\end{eqnarray*} 
are the same, where $\pi$ is the canonical quotient map.
\end{Theorem}
Duflot, Landweber, and Stong point out  \cite{DLS} that a group homomorphism $\rho :V\rightarrow G$ for which $H^*(BV;\mathbb{F}_p)$ becomes a finitely generated $H^*(BG;\mathbb{F}_p)-$module has to be a monomorphism. They pose the question: If $H^*(BG;\mathbb{F}_p)$ is Noetherian and $\phi , \psi :H^*(BG;\mathbb{F}_p)\rightarrow H^*(BV;\mathbb{F}_p)$ are homomorphisms of algebras over the Steenrod algebra $\mathcal{A}^*$ making $H^*(BV;\mathbb{F}_p)$ into a finitely generated $H^*(BG;\mathbb{F}_p)-$module, such that the compositions
\begin{eqnarray*}
H^*(BG;\mathbb{F}_p)\underset{res^*}{\overset{\phi}{\overset{\longrightarrow}{\longrightarrow}}}H^*(BV;\mathbb{F}_p)\overset{\pi }{\rightarrow }H^*(BV;\mathbb{F}_p)/\sqrt{0}
\end{eqnarray*}
are the same then is $\phi =\psi$?
This is answered and more in \cite{Zarati} in the affirmative. One has that.
\begin{Theorem}
If $H^*(BG;\mathbb{F}_p)$ is Noetherian and $\phi , \psi :H^*(BG;\mathbb{F}_p)\rightarrow H^*(BV;\mathbb{F}_p)$ are homomorphisms of algebras over the Steenrod algebra $\mathcal{A}^*$ making $H^*(BV;\mathbb{F}_p)$ into a finitely generated $H^*(BG;\mathbb{F}_p)-$module, such that the compositions
\begin{eqnarray*}
H^*(BG;\mathbb{F}_p)\underset{res^*}{\overset{\phi}{\overset{\longrightarrow}{\longrightarrow}}}H^*(BV;\mathbb{F}_p)\overset{\pi }{\rightarrow }H^*(BV;\mathbb{F}_p)/\sqrt{0}
\end{eqnarray*}
are the same, then $\phi =\psi$.\\
\end{Theorem}
For $H^*(BG;\mathbb{F}_p)$ a Noetherian algebra we know that the minimal prime ideals $p\subset H^*(BG;\mathbb{F}_p)$ are precisely the ideals that arise as the kernels of an induced map $res^*:H^*(BG;\mathbb{F}_p)\rightarrow H^*(BV;\mathbb{F}_p)$, where $V$ ranges over the $p-$elementary abelian subgroups (or in any case their conjugacy classes) $V$ of $G$. This gives the following criterion for nilpotence.
\begin{Theorem}
Suppose that $H^*(BG;\mathbb{F}_p)$ is Noetherian and $u\in H^*(BG;\mathbb{F}_p)$. Then $u$ is nilpotent if and only if for every $p-$elementary abelian subgoup $V$ of $G$ the element $u$ belongs to the kernel of the induced map $Res^*:H^*(BG;\mathbb{F}_p)\rightarrow H^*(BV;\mathbb{F}_p)$. If $G$ has a $p-$Sylow subgroup $S\leq G$ then it suffices to know that $u$ is in the kernel of the induced map $Res^*:H^*(BG;\mathbb{F}_p)\rightarrow H^*(BS;\mathbb{F}_p)$.
\end{Theorem}
\begin{Theorem}
Let $\mathcal{G}$ be a group of finite virtual cohomological dimension with a Sylow $p-$subgroup $S$ and $\mathcal{F}_S(\mathcal{G})$ a saturated fusion system. Then $H^*(B\mathcal{G};\mathbb{F}_p)$ is $F-$isomorphic to the stable elements.
\end{Theorem}
\underline{Proof:} Recall that $\mathcal{G}$ is a group of finite virtual cohomological dimension. It is proved in \cite{Quillen2} that $H^*(B\mathcal{G};\mathbb{F}_p)$ is $F-$isomorphic to the inverse limit of $H^*(-;\mathbb{F}_p)$ over the Quillen category which will be denoted by $H^*(Q)$. Moreover Broto, Levi, and Oliver show in \cite{BLO2} that the ring of stable elements $H^*(\mathcal{F})$ is $F-$isomorphic to the ring $H^*(Q)$. 
 We therefore obtain a commutative diagram:\xymatrix@R=9pt@C=9pt{ {H^*(B\mathcal{G})}\ar[1,1]_{h}\ar[0,1]^{f}&{H^*(\mathcal{F})}\ar[1,0]^{g}\\
&{H^*(Q),}\\
} where $f$ is the natural map in $\mathbb{F}_p-$cohomology of Theorem 2.3, and $g$, $h$ are the natural restriction maps to $H^*(Q)$ respectively. It follows from a diagram chase that the map $f$ is an $F-$isomorphism. Since $ker(f)\subseteq ker (h)$ and $h$ is an $F-$monomorphism it follows that all the elements in $ker(f)$ are nilpotent and therefore $f$ is an $F-$monomorphism as well. It remains to show that $f$ is an $F-$epimorphism, which means that for every $x$ in $H^*(\mathcal{F})$ there is $n\geq 0$ and $y$ in $H^*(B\mathcal{G})$ such that $f(y)=x^n$. Since $h$ is an $F-$epimorphism we have that for every $x$ in $H^*(\mathcal{F})$ there exist $m\geq 0$ and $y'$ in $H^*(B\mathcal{G})$ such that $h(y')=g(x)^m=g(x^m)$. Moreover we have because of the commutativity that $g(f(y')-x^m)=0$. This implies that there is $k\geq 0$ with $(f(y')-x^m)^{k}=0$ and therefore we have $(f(y')-x^m)^{p^k}=0$ and since we are in characteristic $p$ this implies $f(y'^{p^k})=f(y')^{p^k}=x^{mp^k}$. $\Box$
\begin{Theorem}
Let $(S,\mathcal{F},\mathcal{L})$ be a $p-$local finite group and $
\mathcal{G}$
a discrete group such that $S\in Syl_p(\mathcal{G})$ and $\mathcal{F}=\mathcal{F}_S(\mathcal{G})$. Then
there exist a natural map of algebras $H^*(X;A)\overset{q}{\rightarrow} H^*(\mathcal{F};A)$ making $H^*(\mathcal{F};A)$ a module over $H^*(X;A)$ for any system of local coefficients $A$ which we omit in notation.
\end{Theorem}
\underline{Proof:} We show that the restriction map $f^*:H^*(X)\rightarrow H^*(BS)$ factors through the ring of stable elements $H^*(\mathcal{F})\subset H^*(BS)$. For every subgroup $P$ of $S$ we have a map $Res^{X}_P: H^*(X)\rightarrow H^*(BP)$. 
Since we have $\mathcal{F}=\mathcal{F}_{S,f}(X)$ we obtain that for all subgroups $P,  Q\leq S$ and all morphisms $\phi\in Mor_{\mathcal{F}}(P,Q)$ that  the diagram
\xymatrix@R=7pt@C=7pt{ {BP}\ar[0,2]^{B\phi}\ar[1,1]_{Bincl}&{}&{BQ}\ar[1,-1]^{Bincl}\\
&{X}&\\} 
commutes up to homotopy. Therefore the outer triangle in the diagram 
\xymatrix@R=40pt@C=40pt{ &&{H^*(BP)}&{}\\
{H^*(X)}\ar[1,2]_{Res^{\mathcal{G}}_Q}\ar[-1,2]^{Res^{\mathcal{G}}_P}\ar[0,1]^{Res^{\mathcal{G}}_S}&{H^*(BS)}\ar[1,1]^{Res^S_Q}\ar[-1,1]_{Res^S_P}&&\\
&&{H^*(BQ)}\ar[-2,0]^{\phi ^*}\restore&{}\\
} \\
commutes and the map $f^*:H^*(X)\rightarrow H^*(BS)$ factors through the algebra of stable elements giving rise to the map $q:H^*(X)\rightarrow H^*(\mathcal{F})$ which is in particular an epimorphism in the sense of Quillen.$\Box$\\[0.3cm]
The interesting question is under which conditions the map $q$ is an $F$-monomorphism. This is true if $X=B\mathcal{G}$ for a group model $\mathcal{G}$ for $\mathcal{F}$ of finite vcd or if $H^*(B\mathcal{G};\mathbb{F}_p)$ is noetherian as shown above.\\
Moreover we have that if $X=BG\pcom$ for a finite group $G$ or $X=|\mathcal{L}|\pcom$ then the map $q$ is a proper isomorphism as shown by Broto, Levi, Oliver \cite{BLO2}. For $X=B\mathcal{G}\pcom$ for $\mathcal{G}$ a model of the type we have with Leip (which contains the Robinson type and our model with Libman as a special case) the map is an $F$-isomorphism, For $X=B\mathcal{G}\pcom$ with $\mathcal{G}$ a Leary-Stancu model we have a lot less control over the homotopy type since in this case $B\mathcal{G}$ is not $p$-good as was shown in \cite[Proposition 6.9]{anewfam}.\\

We illustrate that the group model does not have to have finite virtual cohomological dimension.
\begin{Proposition}
Let $\mathcal{G}$ be a group model with finite virtual cohomological dimension for a saturated fusion system $\mathcal{F}$ over the finite $p-$group $S$ and $H^*(B\mathcal{G};\mathbb{F}_p)$ is $F-$isomorphic in the sense of Quillen. Let $\mathcal{G}'$ be the group obtained by crossing $\mathcal{G}$ with an infinite number of copies of the integers. Then $\mathcal{G}'$ does not have finite virtual cohomological dimension and is a group model for $\mathcal{F}$. Moreover $H^*(B\mathcal{G}';\mathbb{F}_p)$ is not strongly $F-$isomorphic to the stable elements in the sense of Quillen. 
\end{Proposition}
\underline{Proof:} The fusion systems of $\mathcal{G}$ and $\mathcal{G}'$ are isomorphic. The second part follows from the K\"unneth Formula which is an isomorphism in this case. $\Box$\\[0.3cm]
In certain cases we can say exactly that the cohomology is isomorphic to
all of
 the stable elements.
\begin{Proposition}
Let $X$ be a $p$-complete space, $S$ a finite $p-$group, and $ f: BS \rightarrow X $ a map where
 such that
 $f$ is Sylow,
$f|_{BP}$ is a centric map for each $\mathcal{F}_{S,f}(X)-$centric subgroup $P\leq S$,
and every $\mathcal{F}'_{S,f}(X)-$centric subgroup of $S$ is also $\mathcal{F}_{S,f}(X)-$centric.
 Then $H^*(X;\mathbb{F}_p)$ is isomorphic to the stable elements of the fusion system $\mathcal{F}_{S,f}(X)$ as an unstable algebra over the Steenrod algebra $\mathcal{A}_p$.
\end{Proposition}
\underline{Proof:} The space $X$ is the classifying space of the $p$-local finite group $(S,\mathcal{F}^c_{S,f}(X),\mathcal{L}^c_{S,f}(X))$. $\Box$\\[0.3cm]
The following Proposition not only provides an application to the Kan-Thurston Theorem but also illustrates the limit of fusion to control the kernel of the restriction map in terms of $F-$isomorphism.
\begin{Proposition}
Let $X$ be a $p$-complete space and $ f: BS \rightarrow X $ a map where $BS$ is the classifying space of a finite $p$-group $S$
 such that
 $f$ is Sylow
$f|_{BP}$ is a centric map for each $\mathcal{F}_{S,f}(X)-$centric subgroup $P\leq S$;
and every $\mathcal{F}'_{S,f}(X)-$centric subgroup of $S$ is also $\mathcal{F}_{S,f}(X)-$centric. Let $G$ be the Kan-Thurston group for $\mathbb{C}P^{\infty}$ and $BG$ its classifying space.
 Then $\mathcal{F}_{S,f}(X)=\mathcal{F}_{S,f}(X\times BG)$ and the cohomology of $X\times BG$ is never mapped $F-$monomorphically to the stable elements. 
\end{Proposition}
\underline{Proof:} With the K\"unneth theorem which is an isomorphism in this case since $\mathbb{F}_p$ is a field and since $H^*(BG;\mathbb{F}_p)\cong \mathbb{F}_p[x_2]$ we have $H^*(X\times BG ;\mathbb{F}_p)\cong H^*(\mathcal{F})\otimes _{\mathbb{F}_p} H^*(BG;\mathbb{F}_p)\cong H^*(\mathcal{F})\otimes \mathbb{F}_p[x_2]$. $\Box$\\[0.3cm]
More generally, we have that the kernel of the restriction map can get as arbitrary as possible as long as it is realizable as the cohomology of a space in the category of unstable modules over $\mathcal{A}_p$.
\begin{Theorem}
Let $X$ be a $p$-complete space and $ f: BS \rightarrow X $ a map where $BS$ is the classifying space of a finite $p$-group $S$
 such that\
 $f$ is Sylow
$f|_{BP}$ is a centric map for each $\mathcal{F}_{S,f}(X)-$centric subgroup $P\leq S$;
and every $\mathcal{F}'_{S,f}(X)-$centric subgroup of $S$ is also $\mathcal{F}_{S,f}(X)-$centric. Let $G$ be the Kan-Thurston group for for an arbitrary space $Y$ and $BG$ its classifying space.
 Then $\mathcal{F}_{S,f}(X)=\mathcal{F}_{S,f}(X\vee BG)$ and the cohomology of $X\vee BG$ is isomorphic to $H^*(\mathcal{F})\oplus \overline{H}^*(Y;\mathbb{F}_p)$. 
\end{Theorem}
\underline{Proof:} With the Mayer-Vietoris-Sequence we have for the field $\mathbb{F}_p$ is a field an isomorphism $H^*(BG;\mathbb{F}_p)\cong H^*(Y;\mathbb{F}_p)$ and we have $H^*(X\vee BG ;\mathbb{F}_p)\cong H^*(\mathcal{F})\oplus \overline{H}^*(BG;\mathbb{F}_p)\cong H^*(\mathcal{F})\oplus\overline{H}^*(Y;\mathbb{F}_p)$. $\Box $
\begin{Remark}
Note that there is no assumption made on $\mathcal{F}_{S,f}$ as long as it has a classifying space.
\end{Remark}
\begin{Proposition}
Let $(S,\mathcal{F},\mathcal{L})$ be a $p$-local finite group and $|\mathcal{L}|\pcom$ be its classifying space and $\mathcal{G}$ its Kan-Thurston group.
\end{Proposition}
\underline{Proof:} Assume first that $\mathcal{F}=\mathcal{F}_S(G)$ for $G$ finite. If the order of $G$ is divisible by $p$, then there are non-nilpotent elements of
positive degree in $H^*(G,\mathbb{F}_p)$. Indeed, given any element $x$ of order $p$ in $G$,
there is an element of positive degree cohomology whose restriction to
the cyclic subgroup generated by $x$ is non-nilpotent.In the case that $\mathcal{F}$ is exotic we have that.The Cartan-Eilenberg stable element method still works there,
and one can use Evens norms to produce a stable element.$\Box$
\begin{Proposition}
Let $(S,\mathcal{F},\mathcal{L})$ be a $p$-local finite group and $|\mathcal{L}|\pcom$ be its classifying space and $\mathcal{G}$ its Kan-Thurston group. Then the restriction map $H^*(B\mathcal{G})\rightarrow H^*(Q)$ is neither an $F$-monomorphism nor an $F$-epimorphism.
\end{Proposition}
\underline{Proof:} This follows from the Kan-Thurston construction because all elements are torsion free.$\Box$
\begin{Remark}
The explicit obstruction theory for the existence of a space $Y$ with $H^*(Y;\mathbb{F}_p)\in\mathcal{A}_p$ is described in \cite{Blanc}.
\end{Remark}
\begin{Remark}
It remains open whether the kernel can be isomorphic to a reduced $\mathcal{A}_p-$module.
\end{Remark}
We show that spaces can have the right fusion system without having a nontrivial map from $BS$.
\begin{Proposition}
Let 
$(S,\mathcal{F},\mathcal{L})$ be a $p-$local finite group and $\mathcal{G}$ be the Kan-Thurston group for its either completed or uncompleted classifying space. Then there is no nontrivial map $BS\rightarrow B\mathcal{G}$ and $B\mathcal{G}\cong H^*(\mathcal{F})$. \end{Proposition}
\underline{Proof:} This follows from the Kan-Thurston Construction because the Kan-Thurston group is torsion free.$\Box$
\begin{Theorem}
Let $\mathcal{F}$ be an arbitrary fusion system of a finite $p-$group $S$. Then there is a group model  $\mathcal{G}$ of finite virtual cohomological dimension such that the cohomology of $B\mathcal{G}$ is $F-$isomorphic to $H^*(\mathcal{F})$ for coefficients in a $\mathbb{Z}_{(p)}-$ module $M$.
\end{Theorem}
\underline{Proof:} This follows from the Leary-Stancu group model construction and \cite{LearySchusterYagita}. $\Box$
\begin{Theorem}
Let $\mathcal{F}$ be an arbitrary fusion system of a finite $p-$group $S$ and $Y$ be a space. Then there is a group model  $\mathcal{G}$ of finite virtual cohomological dimension such that the cohomology of $B\mathcal{G}$ is $F-$epimorphic to $H^*(\mathcal{F})$ for coefficients in a $\mathbb{Z}_{(p)}-$ module $M$ and the kernel is isomophic to $\overline{H}^*(Y)/\sqrt{0}$.
\end{Theorem}
\underline{Proof:} This follows from the Leary-Stancu group model construction and \cite{LearySchusterYagita} and the Mayer-Vietoris exact sequence. $\Box$
\begin{Remark}
We provide an example which is not a Quillen group without evoking the Kan-Thurston Theorem. We cite from \cite[Chapter 1]{Henn} that Quillen groups are not closed under colimits.
\end{Remark}
\begin{Proposition}
For the prime $2$ there is an explicit example \cite[Chapter 1]{Henn} that Quillen groups are not closed under colimits.
\end{Proposition}
\underline{Proof:} We can construct an explicit nonnilpotent class in the kernel of the restriction map \cite{AC}, \cite{InoueKono}. $\Box$
\begin{Remark}To the best of our knowledge there is no explicit example known for an odd prime.\end{Remark}
\begin{Theorem}
Let $(S,\mathcal{F},\mathcal{L})$ be a $p-$local finite group and $X$ a space and $f:BS\rightarrow X$ a map such that $\mathcal{F}_{S,f}(X)=\mathcal{F}$ and moreover assume there is a map $g:X\rightarrow |\mathcal{L}|\pcom$ with $f\circ g|_{BS}$ homotopic to $B\delta _S$. Then there is a split short exact sequence of unstable modules over the Steenrod algebra
$0\rightarrow W \rightarrow H^* (X) \overset{\leftarrow}{\rightarrow} H^* (\mathcal{F}) \rightarrow 0,$
with 
$W \cong Ker(f^* )$ in the category of unstable modules over $\mathcal{A}_p$.
\end{Theorem}
\underline{Proof:} The map $q:  H^* (X)\rightarrow  H^* (\mathcal{F})$ was constructed in the previous theorem. The map in the inverse direction $H^*(\mathcal{F})\rightarrow H^*(X)$ is constructed using the isomorphism $H^*(|\mathcal{L}|\pcom)\cong H^*(\mathcal{F})$ as proved in \cite{BLO2}. In fact we have a split short exact sequence $0\rightarrow W \rightarrow H^* (X) \overset{\leftarrow}{\rightarrow} H^* (\mathcal{F}) \rightarrow 0,$
with 
$W \cong Ker(f^* )\in\mathcal{A}_p$ for arbitrary systems of coefficients where for $p-$local coefficients we use \cite{Molinier}.$\Box$\\[0.3cm]
The following proposition illustrates that the cokernel of the restriction map can be arbitrarily big.
\begin{Proposition}
Let $\mathcal{F}$ be the maximal fusion system over a finite $p-$group $S$. Then the fusion system $\mathcal{F}$ is the fusion system of a point $\mathcal{F}=\mathcal{F}_{S,f}(\{*\})$ and the cokernel of the restriction map is maximal.
\end{Proposition}
\underline{Proof:} We have that all injective groups homomorphisms between subgroups are morphisms in the fusion category.$\Box$
\begin{Remark}
It is unknown when the maximal fusion system over a finite $p$-group is saturated.
\end{Remark}
\begin{Remark}
In particular this means that the cokernel does not have to be nilpotent as an algebra. An explicit example is the trivial fusion system of the cyclic group $C_2$ which is maximal.
\end{Remark}
\begin{Proposition}
Let $X$ be a space such that $H^*(X;\mathbb{F}_p)$ is nilpotent and $ f: BS \rightarrow X $ a map where $BS$ is the classifying space of a finite $p$-group $S$
 such that
 $f$ is Sylow
$f|_{BP}$ is a centric map for each $\mathcal{F}_{S,f}(X)-$centric subgroup $P\leq S$;
and every $\mathcal{F}'_{S,f}(X)-$centric subgroup of $S$ is also $\mathcal{F}_{S,f}(X)-$centric. Let $G$ be the Kan-Thurston group for $\mathbb{C}P^{\infty}$ and $BG$ its classifying space.
 Then $\mathcal{F}_{S,f}(X)=\mathcal{F}_{S,f}(X\times BG)$ and the cohomology of $X\times BG$ is never $F-$monomorphic to the stable elements. 
\end{Proposition}
\underline{Proof:} First we have $\mathcal{F}_{S,f}(X)=\mathcal{F}_{S,f}(X\times BG)$ via the construction of the Kan-Thurston group which is torsion free. Second we have that $H^*(X\times BG;\mathbb{F}_p)\cong H^*(X;\mathbb{F}_p)\otimes H^*(\mathbb{C}P^{\infty};\mathbb{F}_p))\cong H^*(\mathcal{F})\otimes \mathbb{F}_p[x_2]$.$\Box$
\begin{Remark}
Similarly we can construct elements of arbitrary degree of nilpotence in the kernel.
\end{Remark}
\begin{Proposition}
Let $X$ be a space such that $H^*(X;\mathbb{F}_p)$ is nilpotent and $ f: BS \rightarrow X $ a map where $BS$ is the classifying space of a finite $p$-group $S$
 such that
 $f$ is Sylow
$f|_{BP}$ is a centric map for each $\mathcal{F}_{S,f}(X)-$centric subgroup $P\leq S$;
and every $\mathcal{F}'_{S,f}(X)-$centric subgroup of $S$ is also $\mathcal{F}_{S,f}(X)-$centric. Let $G$ be the Kan-Thurston group for $\mathbb{C}P^{n}$ and $BG$ its classifying space.
 Then $\mathcal{F}_{S,f}(X)=\mathcal{F}_{S,f}(X\times BG)$ and the kernel of the restriction map in cohomology of $X\times BG$ to the stable elements contains an element of degree of nilpotence $n$. 
\end{Proposition}
\underline{Proof:} First we have $\mathcal{F}_{S,f}(X)=\mathcal{F}_{S,f}(X\times BG)$ via the construction of the Kan-Thurston group which is torsion free. Second we have that $H^*(X\times BG;\mathbb{F}_p)\cong H^*(X;\mathbb{F}_p)\otimes H^*(\mathbb{C}P^{n};\mathbb{F}_p))\cong H^*(\mathcal{F})\otimes \frac{\mathbb{F}_p[x_2]}{(x^n)}$.$\Box$
\begin{Remark}
There is no example known of a group where the restriction homomorphism is a proper $F$-epimorphism.
This has to be an infinite group which has a Sylow $p-$subgroup and a saturated fusion system which has no map to the $p$-competed classifying space restricting up to homotopy to the identity on the image of the embedding of the classifying space of the Sylow $p$-subgroup. \end{Remark}
\begin{Remark}
More general it is unknown which proper subalgebras if any are realizable as the target of the restriction map. For naturality reasons only subalgebras of the stable elements algebra which are closed under Steenrod operations can be the target of the restriction map. But we have.
\end{Remark}
\begin{Proposition}
Let $G$ be a finite group and $R^*\subseteq H^*(BG;\mathbb{F}_2)$ a subalgebra which is closed under the action of the Steenrod algebra. Then there exists a space $X$ and a map $BG\rightarrow X$ such that $Im (f^*)=R^*$.
\end{Proposition}
\underline{Proof:} Choose algebra generators $\{ x_i \}$ for $R^*$ and form the map $f:BG\rightarrow \underset{i\in\mathcal{J}}{\times}K(\mathbb{Z}/2\mathbb{Z}, deg(x_i))$ whose $i-$th component pulls the fundamental class $\iota _{}deg(x_i)$ back to to $x_i$. Then $Im\{f^*H\}=R^*$. $\Box$
\begin{Remark}
The statement holds analogously for odd primes.
\end{Remark}
\begin{Remark}
There is an inclusion of $V=(\mathbb{Z}/2\mathbb{Z})^n\rightarrow \Sigma _{2^n}$ such that the induced map is exactly the Dickson algebra.
\end{Remark}
\begin{Remark}
We quote from \cite{Henn}. The map $q_G$ is not an isomorphism in general. A first counterexample are the cyclic groups $C_{p^k}$, $k>1$, in which case the map is neither injective nor surjective. However, there are interesting classes of groups $G$ for which $q_G$ is an isomorphism or at least a monomorphism. For example, if $p=2$, $q_G$ is an isomorphism if $G=D_{2^n}$, the dihedral group of order $2^n$, or if $G=\Sigma _n$, the symmetric groups on $n$ letters, or $G=O_n(\mathbb{F}_q)$, the symmetric group of the bilinear form $b(x,y)=\sum{x_iy_i}$ on $\mathbb{F}_q^n$ if $\mathbb{F}_q$ is a field of odd order. (The general phenomenon underlying these examples is the following: the property that $q_G$ is an isomorphism is preserved by passing from a group $G$ to its wreath product with $C_2$, and is inherited by any group whose Sylow subgroup has this property). Let $D_8$ be the dihedral group of order $8$. It is known that $H^*(D_8;\mathbb{F}_2)\cong\mathbb{F}_2[x_1,y_1,w_2]/(x_1y_1)$ where the indices give the degree of the elements. The elements $x_1$ and $y_1$ have the property $x_1y_1=0$ reflecting the fact that $D_8$ has two non-conjugate elementary abelian $2-$subgroups of rank $2$. This was generalized by ourselves with Libman \cite{LS} as well as in
\cite{anewfam}.
\end{Remark}
\begin{Proposition}
Let $\mathcal{F}$ be a saturated fusion system over a finite $p$-group $S$ and $\mathcal{G}$ a group model for $\mathcal{F}$ such that we have a short exact sequence: $0\rightarrow K\rightarrow H^*(B\mathcal{G})\rightarrow H^*(\mathcal{F})\rightarrow 0$ of unstable modules, $K\in\mathcal{N}il_k$ for $k$ a fixed integer, $k\geq 1$. For all $s<k$ we have $R_s(H^*(B\mathcal{G}))\cong R_s(H^*(\mathcal{F}))$.
\end{Proposition}
\underline{Proof:}
Let $M$ be an unstable module and $k$ a fixed integer, $k\geq 1$. Since
we have a short exact sequence of unstable modules: $0\rightarrow K\rightarrow H^*(B\mathcal{G})\rightarrow 
H^*(\mathcal{F})
\rightarrow 0$, with $K\in\mathcal{N}il_k$. For all $s<k$ we have that $
R_s(
 H^*(B\mathcal{G})
)\cong R_s(
H^*(\mathcal{F})
)$.
$\Box$
\begin{Remark}
For all the classical group models we have $K\in Nil_1$ and therefore $R_s(
 H^*(B\mathcal{G})
)\cong R_s (
H^*(\mathcal{F})
)$.
\end{Remark}
\begin{Theorem}
Let $X$ be a $1$-connected space such that $H^*(X)$ is of finite dimension in each degree and such that $H^*(X)\in Nil_d$, $d>1$. So $H^*\Omega X\in Nil_{d-1}$ and the unstable module $R_{sd-2}(H^*\Omega X)$.
\end{Theorem}
\underline{Proof:} Let $X$ be a $1$-connected space such that $H^*(X)$ is of finite dimension in each degree and such that $H^*(X)\in Nil_d$, $d>1$. So $H^*\Omega X\in Nil_{d-1}$ and the unstable module $R_{sd-2}(H^*\Omega X)$. $\Box$\\[0.3cm]
\begin{Proposition}
Let $M$ be an unstable module and let $d$ be a fixed integer, $d\geq 1$. Assume there is a short exact sequence $0\rightarrow N\rightarrow M\rightarrow C\rightarrow 0$, with $N\in Nil_d$ and $C\in Nil_{2d}$. Then $M\in Nil_d$, for $d\leq s<2d$ we have $R_s(N)\cong _FR_s(M)$, $R_{sd}(M)\cong _F E$ where $E$ is given via an extension of the form ${0}\rightarrow R_{2d}(N)\rightarrow E\rightarrow L\rightarrow {0}$, where $L$ is a submodule of $R_{2d}(C)$.
\end{Proposition}
\underline{Proof:} This follows from \cite[Proposition 1.9]{Schwartz}. Let $M$ be an unstable module and let $d$ be a fixed integer, $d\geq 1$. Assume there is a short exact sequence $0\rightarrow N\rightarrow M\rightarrow C\rightarrow 0$, with $N\in Nil_d$ and $C\in Nil_{2d}$. Then $M\in Nil_d$, for $d\leq s<2d$ we have $R_s(N)\cong _FR_s(M)$, $R_{sd}(M)\cong _F E$ where $E$ is given via an extension of the form ${0}\rightarrow R_{2d}(N)\rightarrow E\rightarrow L\rightarrow {0}$, where $L$ is a submodule of $R_{2d}(C)$.$\Box$
\subsection{On the cohomology of Park's finite groups realizing fusion systems}
We give a formula for the cohomology of the groups $G=S\wr
\Sigma _{e(X)}$ constructed by Park \cite{Sejong}, \cite{Park2}.
\begin{Proposition}
Let $\mathcal{F}$ be a saturated fusion system over the finite $p$-group $S$ and $G=S\wr
\Sigma _{e(X)}$ Park's model for $\mathcal{F}$ with the characteristic biset $X$. Then there is a second quadrant spectral sequence of $A-$modules converging to $H^*(BG;A)$ for every system of local coefficients $A$ on $BG$.
\end{Proposition}
\underline{Proof:}
We have an extension of groups: $1\rightarrow \underset{e(X)}{\prod}{S}\rightarrow S\wr
\Sigma _{e(X)}\rightarrow \Sigma _{e(X)}\rightarrow 1$ with the natural action of $\Sigma _{e(X)}$ on $\underset{e(X)}{\prod}{S}$ which gives rise to a cohomology Leray-Serre spectral sequence for every system of local coefficients $A$ on $BG$: $\{E_r,d_r\}^{r\geq 0}$ with $E_2^{p,q}\cong H^2(\Sigma _{e(X)};(\underset{e(X)}{\prod}{S};A))$ converging to $H^*(S\wr
\Sigma _{e(X)};A)\cong H^*(BG;A)$ as a spectral sequence of graded modules. 
$\Box$
\begin{Remark}
Note that the cohomology of Park's model is almost never properly isomorphic to the stable elements. This follows from the Martino-Priddy-Conjecture but also from the fact that there is never a proper isomorphism onto the cohomology of a proper subgroup for finite $p-$groups.
\end{Remark}
\begin{Theorem}
Let $\mathcal{F}$ be a fusion system over a finite $p$-group $S$ and $G$ Park's model for $\mathcal{F}$. Then $H^*(BG;A)$ is not $F-$monomorphic to $H^*(\mathcal{F};A)$ for any system of local coefficients $A$.
\end{Theorem}
\underline{Proof:} This follows from the construction of $G=S\wr
\Sigma _{e(X)}$ and the fact that Sylow $p-$subgroups of symmetric groups are iterated wreath products of cyclic groups and the structure of the cohomology of $C_p$. $\Box$
\begin{Theorem}
Let $\mathcal{F}$ be a fusion system over a finite $p$-group $S$ and $G$ Park's model for $\mathcal{F}$. Then $H^*(BG;A)$ is $F-$epimorphic to $H^*(\mathcal{F};A)$ for any system of local coefficients $A$.
\end{Theorem}
\underline{Proof:} This follows from the construction of $G=S\wr
\Sigma _{e(X)}$ and the fact that Sylow $p-$subgroups of symmetric groups are iterated wreath products of cyclic groups and structure of the cohomology of $C_p$. $\Box$

\section{Group models over discrete $p$-toral and profinite groups}
First we extend our previous work on group models to fusion systems over discrete $p$-toral groups.
\subsection{Leary-Stancu model for fusion systems over discrete $p$-toral groups}
Our first result is an analogue of the Leary-Stancu model for fusion systems over discrete $p$-toral groups.
\begin{Theorem}
Let $\mathcal{F}$ be any fusion system over a discrete $p$-toral group $S$ generated by a set of morphisms $\Phi=\{\phi _i:P_i\rightarrow Q_i\}_{i \in I}$. Let $T$ be the free group on the set $I$, and define $\pi_{LS}$ as the quotient of the free product $S*T$ by the relations $t_i^{-1}ut_i=\phi_i(u)$ for all $i$ and for all $u\in P_i$. Then $\pi_{LS}$ is a group model for $\mathcal{F}$.
\end{Theorem}
\underline{Proof:} The proof is totally analogous to \cite[Theorem  2]{Ian+Radu} since there are no finiteness assumptions.$\Box$\\[0.3cm]
The statement of the following theorem depends on the ring of coefficients in the sense that there are no Steenrod operations over cohomology rings with coefficients for example in the rationals $\mathbb{Q}$.
\begin{Theorem} Let $\mathcal{F}$ be a saturated fusion system over the discrete $p-$toral group 
$S$ and $\mathcal{G}$ a group model for $\mathcal{F}$ as constructed above and assume there exists a stable elements theorem for $(S,\mathcal{F},\mathcal{L})$. Then there exist natural maps of algebras over the Steenrod
algebra $q: H^* (B\mathcal{G}) \rightarrow  H^* (\mathcal{F})$ and $s^*: H^*(\mathcal{F}) \rightarrow H^*(B\mathcal{G})$ such that we obtain a split
short exact sequence of unstable modules over the Steenrod algebra
$0\rightarrow W \rightarrow H^* (B\mathcal{G}) \overset{\leftarrow}{\rightarrow} H^* (\mathcal{F}) \rightarrow 0,$
with 
$W \cong Ker(Res^{\mathcal{G}}_S )$ as unstable module.
\end{Theorem}

 Let $\Phi=\{\phi _i:P_i\rightarrow Q_i\}_{i \in I}$ be the set of morphisms used in the construction of $\mathcal{G}$. We will consider them from this point on as morphisms $\phi _i:P_i\rightarrow S$ for all $i\in I$. In \cite{Ian+Radu} the authors show that $B\mathcal{G}=\underset{\mathcal{D}}{hocolim F}$ where $\mathcal{D}$ is the following category:\\
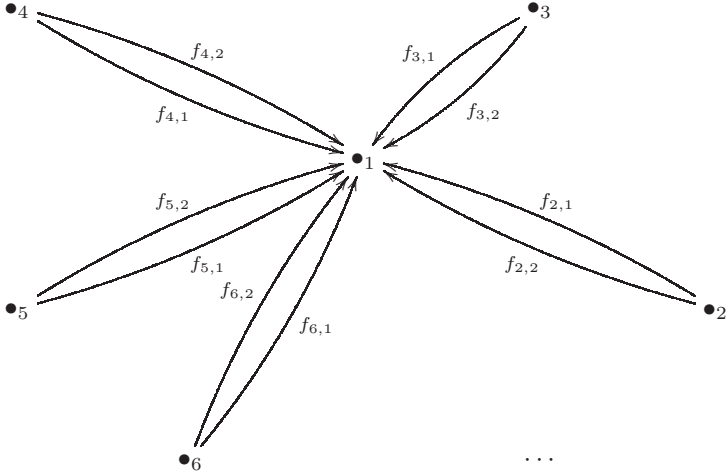
\xymatrix@R=45pt@C=50pt{
{\bullet _4}\ar @/^/[1,2]^{f_{4,2}}\ar @/_/[1,2]_{f_{4,1}}&&&{\bullet _3}\ar @/^/[1,-1]^{f_{3,2}}\ar @/_/[1,-1]_{f_{3,1}}&\\
&&{\bullet _1}&&\\
{\bullet _5}\ar @/_/[-1,2]_{f_{5,1}}\ar @/^/[-1,2]^{f_{5,2}}&&&&{\bullet _2}\ar @/^/[-1,-2]^{f_{2,2}}\ar @/_/[-1,-2]_{f_{2,1}}\\
&{\bullet _{6}}\ar @/_/[-2,1]_{f_{6,1}}\ar @/^/[-2,1]^{f_{6,2}}&&{\cdots}&\\
}\\
and $F$ is a functor to spaces with $F(\bullet _1)=BS$ and $F(\bullet _i)=BP_i$ and $F(f_{i1})=Bincl:BP_i\rightarrow BS$ and $F(f_{i2})=B\phi _i:BP_i\rightarrow BS$ for all $i\in I$.
Due to Alperin's fusion theorem for $p$-local compact groups there exist for all $i\in I$ an index $k(i)$ and a family of $\mathcal{F}-$centric subgroups $P^1_i,\cdots ,P^{k(i)}_i$ and for all $j=1,\cdots,k(i)$ $\psi ^j_i\in Aut_{\mathcal{F}}(P^j_i)$ such that for all $i\in I$, for all $x\in P_i$ we have $\phi _i(x)=\psi ^{k(i)}_i\circ \psi ^{k(i)-1}_i\circ\cdots\circ\psi ^1_i(x)$. Whenever we consider the automorphisms $\psi ^j_i$ as morphisms $\psi ^j_i:P^j_i\rightarrow S$ they will be denoted $\widehat{\psi ^j_i}:P^j_i\rightarrow S$ for all $i\in I$ and $j=1,\cdots ,k(i)$. Note that all the groups $P_i^j$ are $\mathcal{F}-$centric subgroups of $S$ and therefore there is a functor $F_i^{2j}:\mathcal{B}P^j_i\rightarrow\mathcal{L}=incl\circ\delta _{P_i^j}$ for all $i\in I$, $j=1,\cdots ,k(i)$. Therefore we have a map $BP_i^j\rightarrow |\mathcal{L}|$ for all $i\in I$, $j=1,\cdots ,k(i)$. Each object $BP$ of the image of the category $\mathcal{D}$ under $F$ gets mapped to $BAut_{\mathcal{L}}(P)\subset |\mathcal{L}|$. For all $i\in I$ define the family of functors $\{F^j_i:\mathcal{B}P_i\rightarrow \mathcal{L}\}^{1\leq j\leq 2k(i)+1}_{1\leq i\leq n}$ in the following way. For $j=1,\cdots , k(i)$ $F^{2j-1}_i:\mathcal{B}P_i\rightarrow\mathcal{L}$, $x\mapsto \widehat{\psi ^{j-1}_i}\circ\psi ^{j-2}_i\circ\cdots\circ\psi ^1_i(x)$, $\bullet\mapsto S$, $F^{2j}_i:\mathcal{B}P_i\rightarrow\mathcal{L}$, $x\mapsto \widehat{\psi ^{j-1}_i}\circ\psi ^{j-2}_i\circ\cdots\circ\psi ^1_i(x)$, $\bullet\mapsto P_i^j$,$F^{2j+1}_i:\mathcal{B}P_i\rightarrow\mathcal{L}$, $x\mapsto \widehat{\psi ^{j}_i}\circ\psi ^{j-1}_i\circ\cdots\circ\psi ^1_i(x)$, $\bullet\mapsto S$. It follows from the existence of the linking system that we can find lifts of the inclusion $\{\iota _{P_i^j,S}\}$ and lifts of the morphisms $\{\widehat{\psi ^j_i}\}$ which will be denoted $\{\widetilde{\psi ^j_i}\}$ such that for all $i \in I$ and for all $l=1,\cdots 2k(i)+1$ the functors $F_i^{2j}$ and $F_i^{2j+1}$ commute via the following natural transformations respectively. The functor $F_i^{2j}$ commutes to $F^{2j-1}_i$ via $\bullet\mapsto\iota _{P_i^j}$ and $F^{2j}_i$ commutes to $F^{2j+1}_i$ via $\bullet\mapsto\widetilde{\psi ^j_i}$. We
have an algebraic map $q:H^*(B\mathcal{G})\rightarrow H^*(|\mathcal{L}|)$.
The induced diagram on classifying spaces 
commutes up to homotopy for all $i\in I, j=1,2$. We have $s:B\mathcal{G}= \underset{\mathcal{D}}{hocolim^{(1)}}(F)\rightarrow |\mathcal{L}|$ inducing the map $s^*:H^*(|\mathcal{L}|;\mathbb{F}_p)\longrightarrow H^*(B\mathcal{G};\mathbb{F}_p)$. This proves us.
\begin{Proposition}
Let $\mathcal{F}$ be any fusion system over a discrete $p$-toral group $S$ generated by a set of at least two morphisms and $\pi_{LS}$ the group model for $\mathcal{F}$ described above then $B\pi _{LS}$ is $R$-bad.
\end{Proposition}
\underline{Proof:} The space $B\pi _{LS}$ contains a wedge of spheres as a retract which is $R$-bad for any nonzero $
R$. $\Box$
\begin{Remark}
This implies that the kernel $W$ in the previous theorem is nontrivial and the Leary-Stancu Construction is not able to provide a geometric proof of the existence and uniqueness of centric linking systems.
\end{Remark}
\begin{Theorem}
Let $(S,\mathcal{F},\mathcal{L})$ be a $p-$local compact group with strictly positive rank and $\mathcal{G}$ a model of Leary-Stancu type for $\mathcal{F}$. Then $H^*(B\mathcal{G};\mathbb{F}_p)$ is not noetherian.
\end{Theorem}
\underline{Proof:} This follows from \cite[Theorem 4.10.]{gmffs} since the generating set $\Phi$ is not finite in this case. $\Box$
\begin{Theorem}
\label{groupmodelkt}
 Let $K=\{K_1,\dots,K_n\}$ be a generating collection, $\Gamma_K$ the associated graph for some choice of $H$'s and $T$ a generating tree in $\Gamma_K$. Then the amalgam $\pi_{K,T}$ over the graph of groups $T$ is a group with the following properties:
\begin{enumerate}
 \item $S$ is a Sylow $p$-subgroup of $\pi_{K,T}$.
 \item $\mathcal{F}=\mathcal{F}_S(\pi_{K,T})$.
\item $H^*(|\mathcal{L}|,\F_p)$ is a retract of $H^*(B\pi_{K,T})$ in the category of unstable algebras. It is equal to the image of $H^*(\pi_{K,T},\F_p)\rightarrow H^*(S,\F_p)$, the product of any two elements
 in the kernel is zero.
 \item $B\pi_{K,T}$ is $p$-good.
 \item $H^*(B\pi_{K,T})$ is finitely generated.
 \item $|\mathcal{L}|\pcom$ is a stable retract of $(B\pi_{K,T} )\pcom $.
 \item $H^*(B\pi_{K,T})$ is $F-$isomorphic to the stable elements in
 the sense of Quillen.
\end{enumerate}
\end{Theorem}
\underline{Proof:} The proof is based on results of Libman and ourselves \cite{LS} generalizing the case of fusion systems over finite $p$-groups \cite{anewfam}. Consider the category $\mathcal{C}$
consisting of objects $l_0,\dots,l_n$ as well as $l_{i,j}$ for those
$i,j$ such that $K_{i,j}$ is in $T$. Those objects correspond to the
groups mentioned above (i.e. $K_0,\dots K_n$ and the  $K_{i,j}$ in
$T$). In $\mathcal{C}$ let there be a unique morphism from
$l_{i,j}$ to $l_i$ and $l_j$. We have a functor $\gamma$ from this
category into $\hoTop$, which sends $l_i$ to $BK_i$ and $l_{i,j}$ to
$BK_{i,j}$, and the morphisms to those induced by the monomorphisms
$k_{j,i}$ and $k_{i,j}$. We can include each of the classifying spaces into
$|\mathcal{L}|$ (induced by an inclusion of categories
$\mathcal{B}K_i\rightarrow
\mathcal{B}\Aut_\mathcal{L}(P_i)\rightarrow  \mathcal{L}$). We now
want to show that this commutes up to homotopy with the the
morphisms from $\mathcal{C}$, i.e. we want the outer diagram to
commute up to homotopy:
\[\begin{xy}
   \xymatrix{
     & BK_i\ar[rd] &  \\
    BK_{i,j}\ar[ru]\ar[rd]\ar[r] & B\Aut_\mathcal{L}(\langle P_i,P_j \rangle)\ar[r] & |\mathcal{L}| .\\
     & BK_j\ar[ru] &
  }
\label{hocolim_diagram}
 \end{xy}\]

To see this consider the functors
$F_1: \mathcal{B}K_{i,j}\rightarrow\mathcal{B}K_i\rightarrow
\mathcal{L}$ and $F_2:  \mathcal{B}K_{i,j}\rightarrow
\mathcal{B}\Aut_\mathcal{L}(\langle P_i,P_j \rangle \rightarrow
\mathcal{L}$. There is a natural transformation $\eta:
F_1\rightarrow F_2$, which maps the unique object to the morphism
$\widehat{e}=\delta(e)\in \Mor_{\mathcal{L}}(P_i,\langle
P_i,P_j\rangle)$. 
The following diagram commutes by the definition of $k_{j,i}$:
\[\begin{xy}
   \xymatrix{
P_i \ar[r]^{\widehat{e}}\ar[d]_{k_{j,i}(f)=F_1(f)} & \langle P_i,P_j\rangle \ar[d]^{F_2(f)=f}\\
P_i\ar[r]_{\widehat{e}} \ar[r] &\langle P_i,P_j\rangle .
  }
 \end{xy}\]
Thus $\eta$ is a natural transformation. Taking realizations shows that the two triangles in the diagram above are homotopy commutative.
By the universal property, this gives a map
$f: \hocolim \gamma \rightarrow |\mathcal{L}|$. We claim that the composition with
$|\mathcal{L}|\rightarrow |\mathcal{L}|_p^\wedge$ and restriction to
$BS\leq BK_t$ is homotopic to $\theta: BS\rightarrow
|\mathcal{L}|_p^\wedge$ (where we use notation as in \cite{LS}). To see this, consider the following diagram
 \[\begin{xy}
     \xymatrix{
      & & & BK_t\ar[rd] &  \\
     {B}S\ar[rr]^{{B}\delta_S}\ar[rrru]^{{B}\delta_{P_t}}\ar[rrrd]_{{B}\delta_S} & &{B}K_{t,0}\ar[ru]_{{B}k_{0,t}}\ar[rd]^{{B}k_{t,0}}& & |\mathcal{L}|\\
      & & & {B}\Aut_\mathcal{L}(S) \ar[ru] &
    }
  \end{xy}\]
where we set $K_0=\Aut_\mathcal{L}(S)$ for the purpose of this diagram. The morphism $\delta_{P_t}: S\rightarrow \Aut_\mathcal{L}(P_t)$ makes sense because we required $P_t\vartriangleleft S$. The image lies in $K_t$ by another condition on generating trees. Thus the upper left morphism makes sense. The diagram 
\[\begin{xy}
   \xymatrix{
P_t \ar[r]^{\widehat{e}}\ar[d]_{\delta_{P_t}(g)} & S \ar[d]^{\delta_S(g)}\\
P_t\ar[r]_{\widehat{e}} \ar[r] &S
  }
 \end{xy}\]
commutes for all $g\in S$ because $\delta$ is a functor and $\widehat{e}=\delta_{P_t,S}(e)$. Thus $k_{0,t}\circ \delta_S=\delta_{P_t}$, and so the upper triangle strictly commutes.
We have $\langle P_t,S\rangle =S$, so $K_{t,0}$ is a subgroup of $\Aut_\mathcal{L}(S)$, so the restriction $k_{t,0}$ is really the inclusion. Thus the lower triangle also commutes strictly, consisting of morphisms whose restriction to $P_t$ is in $K_t$. 
By the reasoning employed at the beginning of this proof, the right square commutes up to homotopy. Composition with $|\mathcal{L}|\rightarrow|\mathcal{L}|_p^\wedge$ yields the claim.
We have a graph of groups $T$ which we can regard as a category. Note that the functor $B\mathcal{G}: T\rightarrow \hoTop$ is
isomorphic to $\gamma$ discussed above. All the groups involved have
a Sylow $p$-subgroup and the map from a chosen Sylow $p$-subgroup of
$K_{i,j}$ to a Sylow $p$-subgroup of $K_j$ is surjective if $K_{i,j}$ is in $T$. Furthermore, $T$ is a tree and contains a path from $v_0=K_t$ to every other vertex. Thus we can apply \cite[Proposition 3.3]{LS}.
This immediately proves claims 1, 5 and the last part of 4, as well as the fact that  $\hocolim \gamma \simeq B\pi_{K,T}$.
By combining this with the result above, we get a map $f:B\pi_{K,T}\simeq \hocolim \gamma\rightarrow  |\mathcal{L}|_p^\wedge$, whose restriction to $BS$ is homotopic to $\theta$.
We want to use \cite[Theorem 1.1]{LS}, so we need $\mathcal{F}\subset \mathcal{F}_S(\pi_{K,T})$. This is the case, as $\pi_{K,T}$ contains all $K_i$, so $\mathcal{F}_S(\pi)$ contains all $K_i'$. But then $\mathcal{F}_S(\pi_{K,T})$ also contains the fusion system generated by the $K_i'$, which is $\mathcal{F}$, as $K$ is generating. Thus \cite[Theorem 1.1]{LS} yields that the group model $\pi_{K,T}$ also has the properties 2,3 and 4. 
The group $\pi _{K,T}$ is a finite amalgam of finite groups which is generated by elements of $p'$-order and $S$. Let $M$ be the subgroup of $\pi _{K,T}$ generated by all elements of $p'$ -order.
Note that $\pi _{K,T}$ and $S$ surjects on $\pi _{K,T} / M$ and therefore $\pi _{K,T} / M$ is a finite $p$-group. The group $M$ is $p$-perfect since it is generated
by $p'$-elements. Let $X$ be the cover of $B\pi _{K,T}$  with fundamental group $M$. Using the results from
\cite[ VII.3.2]{BK}, we have that $X$ is $p$-good and $X\pcom$
is simply connected. Hence the sequence $X\pcom\rightarrow (B\pi _{K,T})\pcom\rightarrow B(\pi _{K,T} / M)$ is a fibration sequence
and so $(B\pi _{K,T})\pcom$ is $p$-complete by \cite[II.5.2(iv)]{BK}. So $B\pi _{K,T}$ is $p$-good. 
 Recall that
   $H_i\in Syl_p(K_i)$ for all $i=1,...,n$. It follows from \cite[Lemma 2.3.]{BLO1} and
    \cite[Theorem 4.4.(a)]{BLO2} that $H^*(B(K_i))$ is finitely generated over $H^*(|\mathcal{L}|)$ for
    all $i=1,...,n$, and $H^*(|\mathcal{L}|)$ is noetherian as follows from
     \cite[Proposition 1.1. and Theorem 5.8.]{BLO2}. Therefore the Bousfield-Kan spectral sequence
     for $H^*(BG)$ is a spectral sequence of finitely generated $H^*(|\mathcal{L}|)-$modules,
     the $E_2$ term with $E_2^{s,t}= \underset{\mathcal{C}}{lim^s}H^t(F(-);\mathbb{F}_p)$ is concentrated
     in the first two columns and $E_{2}=E_{\infty}$ for placement reasons. Therefore $H^*(B\pi_{K,T})$ is a
      finitely generated module over $H^*(|\mathcal{L}|)$. 
For property 7 recall we have a commutative diagram
\xymatrix@R=9pt@C=9pt{&{\Sigma ^{\infty}BS\pcom}\ar[1,1]^{\Sigma ^{\infty}Bincl\pcom}\ar[1,-1]_{\Sigma ^{\infty}B(\delta _S)\pcom}&{}\\
{\Sigma ^{\infty}|\mathcal{L}|\pcom}&&{\Sigma ^{\infty}(B\pi_{K,T})\pcom.}\ar[0,-2]^{\Sigma ^{\infty}q\pcom}\\
} By Ragnarsson's work  \cite{Ragnarsson} there is a map of spectra $\sigma _{\mathcal{F}}:\Sigma ^{\infty}|\mathcal{L}|\pcom\rightarrow\Sigma ^{\infty}BS \pcom$ such that the composition $ \Sigma ^{\infty}|\mathcal{L}|\pcom\overset{\sigma _{\mathcal{F}}}{\longrightarrow}\Sigma ^{\infty}BS\pcom\overset{\Sigma ^{\infty}B(\delta _S)\pcom}{\longrightarrow}\Sigma ^{\infty}|\mathcal{L}|\pcom$ is the identity. 
Since $\Sigma ^{\infty}B(\delta _S)\pcom\circ \sigma _{\mathcal{F}}
=\Sigma ^{\infty}q\pcom\circ\Sigma ^{\infty}Bincl\pcom\circ \sigma
_{\mathcal{F}}$ we have that the simplicial set $|\mathcal{L}|\pcom$ is a stable retract of
$(B\pi_{K,T})\pcom$. 
Point 7 follows from the definition of $F-$isomorphism and 4.$\Box$\\
\section{Extension of control of fusion to spaces and applications}
In extension and as an application of the Inventiones work of Benson, Grodal, Henke we have that.
\begin{Theorem}[group model analogue to Theorem A, \cite{BGH} ]
Let $\iota :H\leq G$ be group models for a saturated fusion system over a finite $p$-group $S$ of index prime to $p$, and consider the induced map on mod $p$ group cohomology $\iota ^*:H^*(G;\mathbb{F}_p)\rightarrow H^*(H;\mathbb{F}_p)$. If for each $x\in H^*(H;\mathbb{F}_p)$, we have $x^{p^k}\in im(\iota ^*)$ for some $k\geq 0$, then $H$ controls $p$-fusion in $G$.
\end{Theorem}
\underline{Proof:} Let $\iota :H\leq G$ be an inclusion of finite groups of index prime to $p$, $p$ an odd prime, and consider the induced map on mod  $p$ group cohomology $\iota ^*:H^*(G;\mathbb{F}_p)\rightarrow H^*(H;\mathbb{F}_p)$. If for each $x\in H^*(H;\mathbb{F}_p)$, we have $x^{p^k}\in im(\iota ^*)$ for some $k\geq 0$, then $H$ controls $p$-fusion in $G$.$\Box$
\begin{Theorem}[group model analogue to Theorem B, \cite{BGH}]
Let $\mathcal{G}$ and $\mathcal{F}$ be group models for a saturated fusion system on the same finite $p-$group $S$. Suppose that
$Hom_{\mathcal{G}}(A,B)=Hom_{\mathcal{F}}(A,B)$ for all $A,B\leq S$ with $A,B$ elementary abelian if $p$ is odd, and abelian of exponent at most $4$ if $p=2$. Then $\mathcal{G}=\mathcal{F}$.
\end{Theorem}
\underline{Proof:} Let $\mathcal{G}\leq\mathcal{F}$ be two group models for a saturated fusion systems on the same finite $p-$group $S$. Suppose that
$Hom_{\mathcal{G}}(A,B)=Hom_{\mathcal{F}}(A,B)$ for all $A,B\leq S$ with $A,B$ elementary abelian if $p$ is odd, and abelian of exponent at most $4$ if $p=2$. Then $\mathcal{G}=\mathcal{F}$.$\Box$
\section{Group cohomology and the centralizer spectral sequence}

We compute the cohomology of the classifying space of a group model via the spectral sequence. 
Jackowski and McClure \cite{JM}
show that we
 have a spectral sequence $hocolim$ $BC_G(E)\rightarrow BG$ and a map $hocolim$ $EG \underset{G}{\times}(G\underset{C_G(E)}{\times}X^E)\rightarrow EG\underset{G}{\times }X_{p,S}$ and a map $EG\underset{G}{\times}X_{p,S}\rightarrow hocolim EG \underset{G}{\times}(G\underset{C_G(E)}{\times}X^E)$. We therefore obtain a second quadrant spectral sequence in homology.
\begin{Theorem}
Let $G$ be a group model for $\mathcal{F}$. Then the natural map
\begin{eqnarray}
\pi :hocolim _{\mathcal{A}^*_p(G)^{op}}EG\times _G G/C_G(E)\rightarrow BG
\end{eqnarray}
induces an isomorphism in homology with coefficients in the $p$-local integers $\mathbb{Z}_{(p)}$.
\end{Theorem} 
\underline{Proof:} This follows from \cite{JM}. $\Box$
\section{A functor to the category of groups}
We can extend our results from \cite{gmffs} as follows. Let $p$ be a prime. Define the category
$FUSION
(p)$. The objects of this category are fusion systems over discrete 
$p$
-toral groups
and its morphisms are the morphisms between the respective fusion systems. Let
$GROUP_
{Syl_
p}$
be the full subcategory of the
category of groups where the objects are groups which have a Sylow
$p$
-subgroup. Define the functor
$F
:
FUSION
(
p
)\rightarrow
GROUP_
{Syl_
p}$, as constructed in
\cite[Corollary 4]{Ian+Radu}. Let
$\mathcal{F}$
be an object of
$FUSION
(
p
)$, i.e. a fusion system over a discrete
$p$
-toral group
$S$.
Then the functor
$F$
takes
$\mathcal{F}$
to the group
$S*
F
(
Mor
(
\mathcal{F}
))/ < \phi
u
\phi
-
1
=
\phi(
u
)\text{ }
\forall
\text{ }
\phi
\in
Mor
(
\mathcal{F}
)
,
\phi
:
P\rightarrow
Q
,
u
\in
P
>$, where
$Mor
(
\mathcal{F}
)$
is the set of all morphisms in
$\mathcal{F}$
and
$F
(
Mor
(
\mathcal{F}
))$
is the free group on the morphism set
$Mor
(
\mathcal{F}
)$. Let
$\mathcal{F}$
and
$\mathcal{F}'$
be fusion systems over the discrete 
$p$
-toral groups
$S$
and
$S'$
respectively. Let
$(
\alpha
,\phi )
:
\mathcal{F}
\rightarrow
\mathcal{F}'$
be a morphism of fusion systems
between them. Define
$F
((
\alpha
,\phi))
:
F
(
\mathcal{F}
)
\rightarrow
F
(
\mathcal{F}'
)$
by
$s
\mapsto
\alpha
(
s
)$
and
$\phi
\mapsto
\Phi(\phi )$. 
The functor $F$
is a left inverse to the canonical functor in $p$-local homotopy theory 
\cite[Corollary 4]{Ian+Radu}. Recall that there is no left adjoint to the canonical functor in fusion theory \cite[Remark 3.2.]{gmffs}. The fact that the canonical functor extends to fusion systems over discrete $p-$toral  and profinite groups implies that there is no left adjoint either to the functor which assigns to a group with a discrete $p$-toral or profinite Sylow $p$-subgroup its fusion system. $\Box$


{\noindent
Dr.\ Nora Seeliger PhD, Department of Mathematics and Statistics, Room B9, Bailrigg Campus, Lancaster University, LA1 4YF, Email: s.nora@lancaster.ac.uk.}

\end{document}